\documentclass[twocolumn,english]{IEEEtran}
\usepackage[T1]{fontenc}
\usepackage{babel}
\usepackage{float}
\usepackage{units}
\usepackage{amsmath}
\usepackage{amssymb}
\usepackage{graphicx}
\usepackage{color}

\makeatletter

\floatstyle{ruled}
\newfloat{algorithm}{tbp}{loa}
\providecommand{\algorithmname}{Algorithm}
\floatname{algorithm}{\protect\algorithmname}

\newtheorem{theorem}{Theorem}
\newtheorem{claim}{Claim}

\newtheorem{lemma}{Lemma}
\newtheorem{problem}{Problem}

\makeatother

\begin{document}

\title{Towards Analog Memristive Controllers}

\author{Gourav~Saha,~Ramkrishna~Pasumarthy,~and Prathamesh Khatavkar%
\thanks{The authors are with the Department of Electrical Engineering, Indian
Institute of Technology, Madras, India. e-mail: \{ee13s005, ramkrishna, ee13s007\}@ee.iitm.ac.in%
}}
\maketitle
\begin{abstract}
Memristors, initially introduced in the 1970s, have received increased
attention upon successful synthesis in 2008. Considerable work has
been done on its modeling and applications in specific areas, however,
very little is known on the potential of memristors for control applications.
Being nanoscopic variable resistors, one can consider its use in making
variable gain amplifiers which in turn can implement gain scheduled
control algorithms. The main contribution of this paper is the development
of a \textit{generic} memristive analog gain control framework and
theoretic foundation of a gain scheduled robust-adaptive control strategy
which can be implemented using this framework. Analog memristive controllers
may find applications in control of large array of miniaturized devices
where robust and adaptive control is needed due to parameter uncertainty
and ageing issues.\end{abstract}

\section{Introduction\label{sec:Intro}}

\IEEEPARstart{M}{emristor} \cite{CHUA1971}, considered as the
fourth basic circuit element, remained dormant for four decades until
the accidental discovery of memristance in nanoscopic crossbar arrays
by a group of HP researchers \cite{STRUKOV2008}. Memristor, an acronym
for memory-resistor, has the capability of memorizing its history
even after it is powered off. This property makes it a desirable candidate
for designing high density non-volatile memory \cite{MEMORY2011}.
However, optimal design of such hardware architectures will require
accurate knowledge of the nonlinear memristor dynamics. Hence, considerable
effort has been channeled to mathematically model memristor dynamics
(see \cite{STRUKOV2008,JOGELKAR2009,BCM2012,TEAM2013}). The memorizing
ability of memristor has lead researchers to think about its possible
use in neuromorphic engineering. Memristors can be used to make dense
neural synapses \cite{SynapseMemristor} which may find applications
in neural networks \cite{memneural}, character recognition \cite{CharacterRecognition},
emulating evolutionary learning (like that of Amoeba \cite{amoebamemristivelearning}).
Other interesting application of memristor may include its use in
generating higher frequency harmonics which can be used in nonlinear
optics \cite{harmonic_generation} and its use in making programmable
analog circuits \cite{progAnalog1,progAnalog2}.

Memristor is slowly attracting the attention of the control community.
Two broad areas have received attention: \textbf{1)} Control of Memristive
Systems. This include works reported in \cite{jeltsema2012port,jeltsema2010memristive}
give detailed insight into modelling and control of memristive systems
in Port-Hamiltonian framework while State-of-the-Art work may include
\cite{stabilityMemNeural} which studies global stability of Memristive
Neural Networks. \textbf{2)} Control using Memristive Systems. The
very first work in this genre was reported in \cite{DELGADO2010}
where the author derived the describing function of a memristor which
can be used to study the existence of undesirable limit cycles (i.e.
sustained oscillations) in a closed loop system consisting of memristive
controller and linear plant. Another line of work may include \cite{MEMIDAPBC2012,cdc_memristor}
which uses memristor as a passive element to inject variable damping
thus ensuring better transient response. In this paper we lay the
groundwork to use memristor as an analog gain control (AGC) element
for robust-adaptive control of miniaturized systems.

\textbf{Why "Analog Memristive Controller"?: }Several applications
needs controlling an array of miniaturized devices. Such devices demands
Robust-Adaptive Control due to parameter uncertainty (caused by design
errors) and time varying nature (caused by ageing effect). Robust-Adaptive
Control algorithms found in literature are so complex that they require
micrcontroller for implementation. This poses scalability and integration
issues \cite{shapiro2011feedback} (page 190) because microcontroller
in itself is a complicated device. The motive of this work is two
folded: 1) Invoke a thought provoking question: ``Can modern control
laws (like Robust-Adaptive control) be implemented using analog circuits%
\footnote{A microcontroller (or a few microcontrollers) may be used to control
an array of miniaturized devices by using Time Multiplexing. In time
multiplexing a microcontroller controls each element of the array
in a cyclic fashion. In such a scenario a microcontroller will face
a huge computational load dependent on the size of the arrray. Upcoming
ideas like ``event-based control'' promises to reduce the computational
load by allowing aperodic sampling. The motive of this work is not
to challenge an existing idea but to propose an alternative one.%
}?''. 2) Suggest memristor as an AGC element%
\footnote{CMOS based hybrid circuits like that proposed in \cite{tunableresistor}
can also act like variable gain control element. However memristors
are much smaller (found below $10\, nm$) than such hybrid circuits
and hence ensures better scalability. Gain-Span of memristor is also
more than CMOS based hybrid circuits.%
} for implementing Robust-Adaptive Control.

\begin{figure}[t]
\begin{centering}
\includegraphics[width=8cm,height=4.4cm]{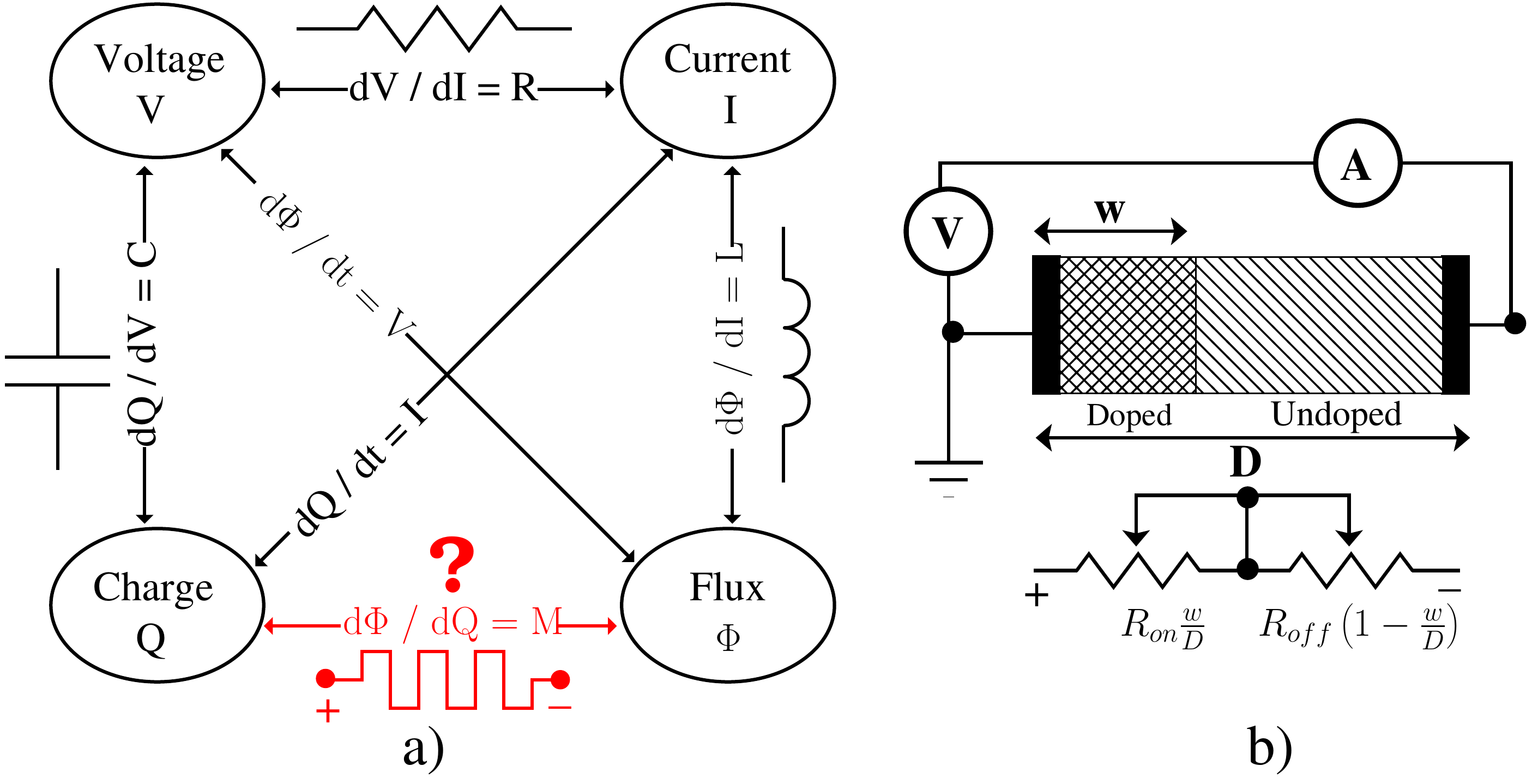}
\par\end{centering}

\caption{a) Memristor as a missing link of Electromagnetic Theory. Adapted from \cite{JOGELKAR2009}. b) Memristor as a series resistor formed by doped and undoped regions of $TiO_{2}$. Adapted from \cite{STRUKOV2008}.}

\vspace{-1.5em}
\end{figure}

The paper is organized as follows. We begin our study by gaining basic
understanding of memristor in Section~\ref{sec:Prelim}. A generic
gain control architecture using memristor is proposed in Section~\ref{sec:memgain}.
Section~\ref{sec:cntrl} discusses a robust-adaptive control strategy
which can be implemented in an analog framework. Section~\ref{sec:Example}
deals with designing an analog controller for a miniaturized setup
by using results from Section~\ref{sec:memgain} and \ref{sec:cntrl}.

\section{Memristor Preliminaries\label{sec:Prelim}}

Chua \cite{CHUA1971} postulated the existence of memristor as a missing
link in the formulation of electromagnetism. Fig. 1a) gives an intuitive
explanation of the same. As evident from Fig. 1a), such an element
will link charge $Q$ and flux $\Phi$, i.e. $\Phi=f\left(Q\right)$.
Differentiating this relation using chain rule and applying Lenz's
Law yields

\vspace{-1.2em}

\begin{equation}
V=M\left(Q\right)I
\end{equation}

\vspace{-0.4em}

\noindent suggesting that the element will act like a charge controlled
variable resistor with $M\left(Q\right)=\frac{df}{dQ}$ as the variable
resistance. This device has non-volatile memory (\cite{CHUA1971,JOGELKAR2009})
and are hence called memristors. Memristive systems \cite{KANG1976}
are generalization of memristor with state space representation,

\begin{figure}[t]
\begin{centering}
\includegraphics[width=0.8\columnwidth,height=0.16\paperheight]{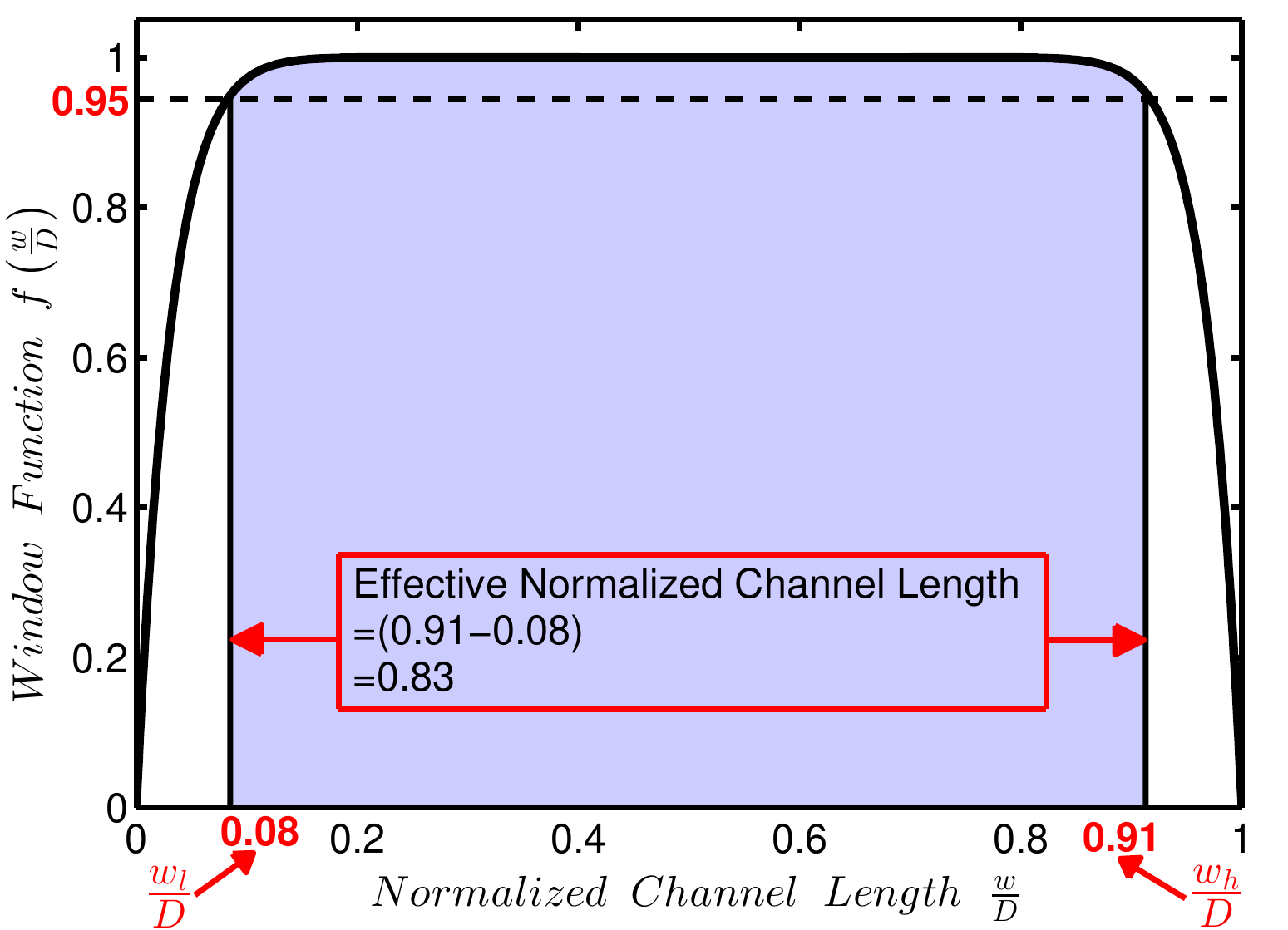}
\par\end{centering}

\caption{Plot of jogelkar window function with $p=8$ for HP memristor. For
a tolerance of $5\%$ in $f\left(\frac{w}{D}\right)$the safe zone
of operation is $\frac{w}{D}\in\left[0.08\,,\,0.91\right]$. }

\vspace{-1.2em}
\end{figure}

\vspace{-1.0em}

\begin{equation}
\dot{W}=F\left(W\,,\, I\right);\; V=R\left(W\,,\, I\right)I
\end{equation}

where, $W$ are the internal state variables, $I$ is the input current,
$V$ is the output voltage, $R\left(W\,,\, I\right)$ is the memristance.
Memristor reported by HP Labs in \cite{STRUKOV2008} is essentially
a memristive system with the following state space representation,

\vspace{-1.0em}

\begin{equation}
\begin{array}{l}
\dot{w}=\mu\frac{R_{on}}{D}f\left(\frac{w}{D}\right)I;\; V=\left[R_{on}\frac{w}{D}+R_{off}\left(1-\frac{w}{D}\right)\right]I\end{array}
\end{equation}

It consists of two layers: a undoped region of oxygen-rich $TiO_{2}$
and doped region of oxygen-deficient $TiO_{2-x}$. Doped region has
channel length $w$ and low resistance while undoped region has channel
length $D-w$ and high resistance. These regions form a series resistance
as shown in Fig. 1b. $R_{on}$ and $R_{off}$ are the effective resistance
when $w=D$ and $w=0$ respectively with $R_{off}\gg R_{on}$. $\mu$
is the ion mobility. When external bias is applied, the boundary between
the two regions drift. This drift is slower near the boundaries, i.e.
$\dot{w}\rightarrow0$ as $w\rightarrow0$ or $w\rightarrow D$. This
nature is captured in \cite{JOGELKAR2009,BCM2012} using window function
$f\left(\frac{w}{D}\right)$. A plot of Joglekar window function $\ensuremath{f\left(\frac{w}{D}\right)=1-\left(\frac{2w}{D}-1\right)^{p}}$
is shown in Fig. 2. A close investigation of various proposed models
\cite{STRUKOV2008,JOGELKAR2009,BCM2012,TEAM2013} reveals two important
facts: 1) As shown in Fig. 2, $f\left(\frac{w}{D}\right)$ is approximately
$1$, except near the boundaries. 2) Boundary dynamics of memristor
is highly non-linear and still a matter of debate. Hence the region
$w\in\left[w_{l},w_{h}\right]$, $0<w_{l}<w_{h}<D$, where $f\left(\frac{w}{D}\right)\thickapprox1$
is the \textit{safe zone} in which memristor dynamics can be approximated
as

\vspace{-0.6em}

\begin{equation}
\dot{Q}_{M}=I\,;\quad V=\left[R_{off}^{S}-\alpha^{S}Q_{M}\right]I
\end{equation}

\noindent where, $\alpha^{S}=\frac{\left(R_{off}^{S}-R_{on}^{S}\right)}{Q_{M}^{S}}$,
$R_{off}^{S}=R_{off}-\left(R_{off}-R_{on}\right)\frac{w_{l}}{D}$,
$R_{on}^{S}=R_{off}-\left(R_{off}-R_{on}\right)\frac{w_{h}}{D}$,
$D^{S}=w_{h}-w_{l}$, $Q_{M}^{S}=\frac{D^{S}D}{\mu R_{on}}$. Superscript
``S'' means ``safe''. In equation (4) \textit{we define} $Q_{M}$
such that $Q_{M}=0$ when $w=w_{l}$. Then $Q_{M}=Q_{M}^{S}$ when
$w=w_{h}$. Hence equation (4) is valid when $Q_{M}\in\left[0\,,\, Q_{M}^{S}\right]$.
From now on, the following conventions will be used:

\noindent \begin{enumerate}
\item Memristor would mean a HP memristor operating in safe zone. Memristor dynamics is governed by equation (4).
\item The schematic symbol of the memristor shown in Fig. 1a will be used. Conventionally, resistance of memristor decreases as the current enters from the port marked "+".
\item HP memristor parameters: $R_{on}=100\Omega$, $R_{off}=16k\Omega$, $D=10nm$, $\mu=10^{-14}\frac{m^{2}}{sV}$. We model the memristor using Joglekar Window Function with $p=8$. The safe zone is marked by $w_{l}=0.08D$ and $w_{h}=0.91D$ (see Fig. 2) in which $f\left(\frac{w}{D}\right)\geq0.95$. This gives: $R^{S}_{off}=15k\Omega$, $R^{S}_{on}=1.5k\Omega$, $Q^{S}_{M}=83 \mu C$, $\alpha^{S}=1.6 \times 10^{8} \frac{\Omega}{C}$. These parameters will used for design purposes.
\end{enumerate}

\section{Analog Gain Control Framework\label{sec:memgain}}

In this section we design a AGC circuit whose input-output relation
is governed by the following equation,

\vspace{-0.5em}

\begin{equation}
\dot{K}=\alpha_{k}V_{C}\left(t\right)\,;\; V_{u}\left(t\right)=KV_{e}\left(t\right)
\end{equation}

We assume that $K>0$. This is basically a variable gain proportional
controller with output voltage $V_{u}$, error voltage $V_{e}$ and
variable gain $K$. $K$ is controlled by control voltage $V_{C}$.
$\alpha_{k}$ determines the sensitivity of $V_{C}$ on $K$. An analog
circuit following equation (5) is \textit{generic} in the sense that
it can be used to implement any gain scheduled control algorithm.
We assume $V_{C}$ and $V_{e}$ are band limited, i.e. - the maximum
frequency component of $V_{C}$ and $V_{e}$ are $\omega_{C}^{M}$
and $\omega_{e}^{M}$ respectively. Knowledge of $\omega_{C}^{M}$
and $\omega_{e}^{M}$ is assumed.

The proposed circuit is shown in Fig. 3. We assume the availability
of positive and negative power supply%
\footnote{$V_{DD}$ and $-V_{DD}$ are the highest and the lowest potential
available in the circuit. From control standpoint, it imposes bounds
on control input $V_{u}$.%
} $V_{DD}$ and $-V_{DD}$. All op-amps are powered by $V_{DD}$ and
$-V_{DD}$.

\noindent \vspace{-1.0em}
\begin{claim}
The proposed circuit shown in Fig. 3. is an \textit{approximate} analog realization of equation (5) if:
\begin{enumerate}
\item Electrical components in Fig. 3. are ideal. Also for MOSFET's, threshold voltage $V_{th}\approx0$.
\item $\omega_{m}=1000\max\left\{ \left(\omega_{C}^{M}+\omega_{e}^{M}\right),2\omega_{C}^{M}\right\}$
\item $R_{f}C_{f}=\frac{100}{\omega_{m}}\,;\; R_{e}C_{e}=\frac{1}{2\left(\omega_{C}^{M}+\omega_{e}^{M}\right)}$
\end{enumerate}
\end{claim}

We understand the working of the proposed circuit by studying its
four distinct blocks and in the process prove the above claim. The
output response of each block for a given input, $V_{e}$ and $V_{C}$,
is shown in Fig. 4. It should be noted that the tuning rules proposed
in the claim is only one such set of values which will make the circuit
follow equation (5).

\textit{Remark 1:} Substrate of all NMOS and PMOS are connected to
$-V_{DD}$ and $V_{DD}$ respectively. In $ON$ state%
\footnote{With a slight abuse of terminology, a NMOS and a PMOS is said to be
in $ON$ state when its gate voltage is $V_{DD}$ and $-V_{DD}$ respectively.%
}, voltages between $-V_{DD}$ to $\left(V_{DD}-V_{th}\right)$ will
pass through%
\footnote{A voltage is said to pass through a MOSFET if the exact voltage applied
at its source(drain) terminal appears at its drain(source) terminal.%
} NMOS and voltages between $-\left(V_{DD}-V_{th}\right)$ to $V_{DD}$
will pass through PMOS. If $V_{th}\approx0$, any voltage between
$-V_{DD}$ to $V_{DD}$ will pass through NMOS and PMOS when they
are in $ON$ state.

\begin{figure}[t]
\begin{centering}
\includegraphics[scale=0.5]{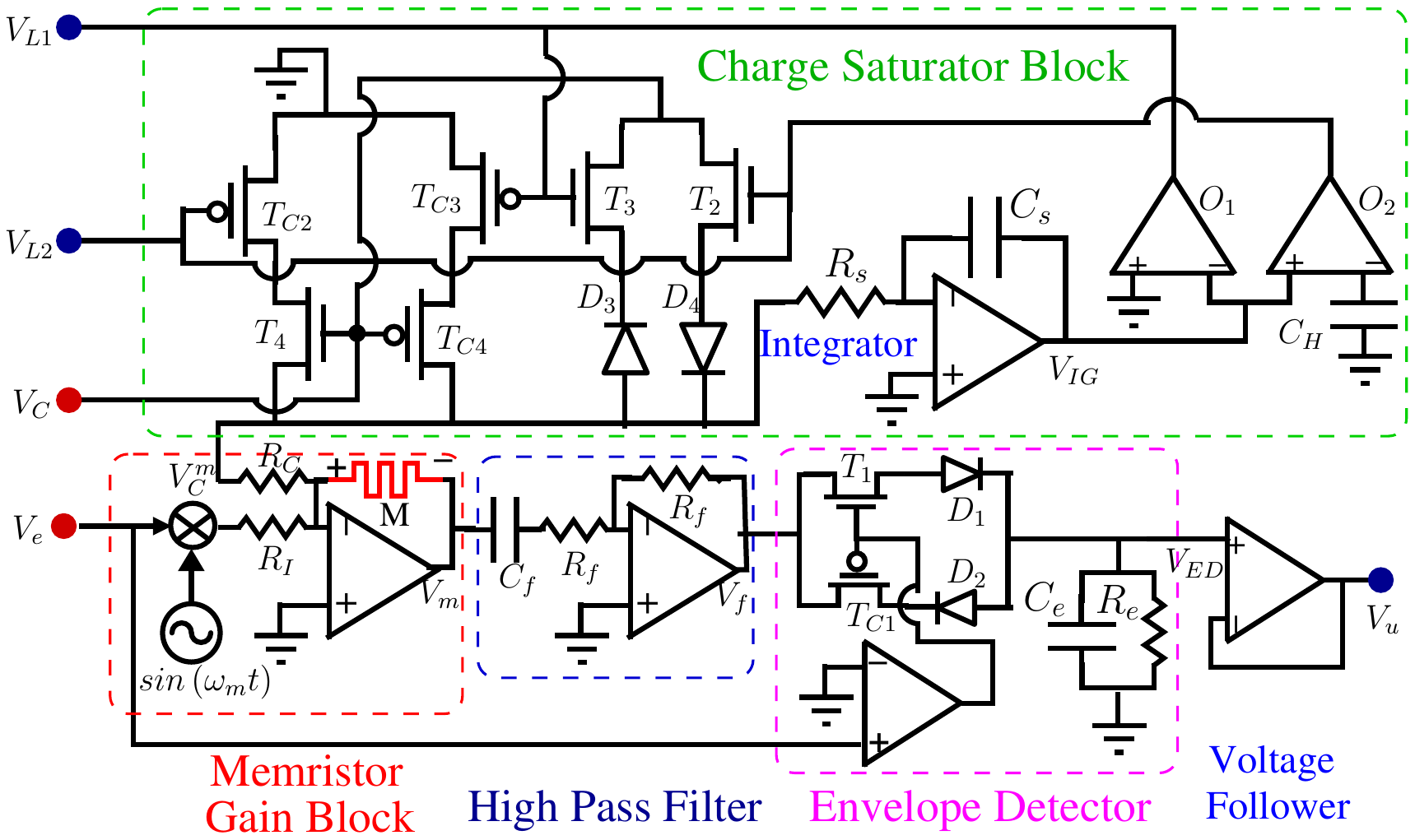}
\par\end{centering}

\caption{Memristive AGC Circuit. $V_{e}$ and $V_{C}$ are the inputs. $V_{u}$, $V_{L1}$ and $V_{L2}$ are the outputs. $V_{L1}$ and $V_{L2}$ are zone indicating voltages.}

\vspace{-1.2em}
\end{figure}

\subsection{Memristor Gain Block\label{sub:Memristor-Gain-Block}}

The key idea of this block has been adapted from \cite{memgain}.
$V_{C}^{m}$ and $V_{e}$ are inputs to this block while $V_{m}$
is the output. For now we assume $V_{C}^{m}=V_{C}$ (details discussed
in \ref{sub:Charge-Saturator-Block}). Current $I_{m}$ through memristor
is

\vspace{-0.4em}

\[
I_{m}\left(t\right)=\underset{I_{noise}}{\underbrace{\frac{V_{e}\left(t\right)\sin\left(\omega_{m}t\right)}{R_{I}}}}+\underset{I_{cntrl}}{\underbrace{\frac{V_{C}\left(t\right)}{R_{C}}}}
\]

\vspace{-0.2em}

From equation (4), resistance of the memristor is given by $M\left(Q_{M}\right)=R_{off}^{S}-\alpha^{S}Q_{M}$.
Differentiating this relation we get, $\dot{M}=-\alpha^{S}I_{m}\left(t\right)$.
Hence voltage $V_{m}$ is given by

\vspace{-0.6em}

\begin{equation}
\dot{M}=-\alpha^{S}I_{m}\left(t\right)\,;\; V_{m}\left(t\right)=-I_{m}\left(t\right)M
\end{equation}

Note that $\omega_{m}=1000\max\left\{ \left(\omega_{C}^{M}+\omega_{e}^{M}\right),\,2\omega_{C}^{M}\right\} >\omega_{e}^{M}$.
In such a case the minimum component frequency of $I_{noise}=\frac{V_{e}\left(t\right)\sin\left(\omega_{m}t\right)}{R_{I}}$
is $\omega_{m}-\omega_{e}^{M}$. Now

\vspace{-1.4em}

\begin{eqnarray*}
\omega_{m}-\omega_{e}^{M} & = & 1000\max\left\{ \left(\omega_{C}^{M}+\omega_{e}^{M}\right),\,2\omega_{C}^{M}\right\} -\omega_{e}^{M}\\
 & \geq & 1000\left(\omega_{C}^{M}+\omega_{e}^{M}\right)-\omega_{e}^{M}\\
 & = & 1000\omega_{C}^{M}+999\omega_{e}^{M}\gg\omega_{C}^{M}
\end{eqnarray*}

This implies that the lowest component frequency of $I_{noise}$ is
much greater than the highest component frequency of $I_{cntrl}$.
Also note that $\dot{M}=-\alpha^{S}I_{m}\left(t\right)$ in an integrator
with input $I_{m}\left(t\right)$ and output $M\left(t\right)$. As
integrator is a low pass filter, the effect of high frequency component
$I_{noise}$ on $M$ is negligible compared to $I_{cntrl}$. Hence
equation (6) can be modified as

\vspace{-1.4em}

\begin{equation}
\dot{M}\thickapprox-\frac{\alpha^{S}}{R_{C}}V_{C}\left(t\right)\,;\; V_{m}\left(t\right)=V_{m}^{1}+V_{m}^{2}\quad\text{{where,}}
\end{equation}

\noindent $V_{m}^{1}=-\frac{MV_{e}\left(t\right)}{R_{I}}\sin\left(\omega_{m}t\right)$
and $V_{m}^{2}=-\frac{MV_{C}\left(t\right)}{R_{C}}$. Note $V_{m}^{1}$
is the modulated form of the desired output with gain $K=\frac{M}{R_{I}}$.

\textit{Remark 2: }$M\left(t\right)$ is the variable gain. According
to equation (5), $V_{e}$ should not effect $M\left(t\right)$. Without
modulating $V_{e}$, the effect of $V_{e}\left(t\right)$ on $M\left(t\right)$
would not have been negligible.

\subsection{High Pass Filter (HPF)\label{sub:High-Pass-Filter}}

The role of HPF is to offer negligible attenuation to $V_{m}^{1}$
and high attenuation to $V_{m}^{2}$ thereby ensuring that the Envelope
Detector can recover the desired output. 

Note that $M\left(t\right)$ is basically the integral of $V_{C}\left(t\right)$.
Since integration is a linear operation it does not do frequency translation.
Hence the component frequencies of $M\left(t\right)$ and $V_{C}\left(t\right)$
are same. Let $\omega_{1}^{m}$ and $\omega_{2}^{M}$ denote the minimum
component frequency of $V_{m}^{1}$ and maximum component frequency
of $V_{m}^{2}$ respectively. Now

\vspace{-1.4em}

\begin{eqnarray}
\omega_{1}^{m} & = & \omega_{m}-\text{{Maximum\:\ Frequency\:\ Component\:\ of\:}}MV_{e}\nonumber \\
 & = & \omega_{m}-\text{{Maximum\:\ Frequency\:\ Component\:\ of\:}}V_{C}V_{e}\nonumber \\
 & = & \omega_{m}-\left(\omega_{C}^{M}+\omega_{e}^{M}\right)\\
\omega_{2}^{M} & = & \text{{Maximum\:\ Frequency\:\ Component\:\ of\:}}MV_{C}\nonumber \\
 & = & \text{{Maximum\:\ Frequency\:\ Component\:\ of\:}}\left(V_{C}\right)^{2}\nonumber \\
 & = & 2\omega_{C}^{M}
\end{eqnarray}

The attenuation offered by the HPF is given by $-10\log\left(1+\left(\omega R_{f}C_{f}\right)^{-2}\right)\, dB$.
We want to study the attenuation characteristics at two frequencies:

1) At $\omega=\omega_{1}^{m}$: At this frequency

\vspace{-1.2em}

\begin{eqnarray}
 &  & \omega_{1}^{m}R_{f}C_{f}\nonumber \\
 &  & =\frac{100\left(\omega_{m}-\left(\omega_{C}^{M}+\omega_{e}^{M}\right)\right)}{\omega_{m}}=100\left(1-\frac{\omega_{C}^{M}+\omega_{e}^{M}}{\omega_{m}}\right)\nonumber \\
 &  & \geq100\left(1-\frac{1}{1000}\right)\approx100
\end{eqnarray}

Inequality (10) is possibe because $\omega_{m}=1000\max\left\{ \left(\omega_{C}^{M}+\omega_{e}^{M}\right),\,2\omega_{C}^{M}\right\} $
implying $\frac{\omega_{m}}{\omega_{C}^{M}+\omega_{e}^{M}}\geq1000$.
For $\omega_{1}^{m}R_{f}C_{f}\geq100$ attenuation is approximately
$0\, dB$. As the HPF offers almost no attenuation to the minimum
component frequency of $V_{m}^{1}$, it will offer no attenuation
to the higher component frequency of $V_{m}^{1}$ as well. Hence $V_{m}^{1}$
suffers negligible attenuation.

2) At $\omega=\omega_{2}^{M}$: At this frequency

\vspace{-0.4em}

\begin{equation}
\omega_{2}^{M}R_{f}C_{f}=100\frac{2\omega_{C}^{M}}{\omega_{m}}\leq\frac{100}{1000}=0.1
\end{equation}

Inequality (11) is possibe because $\omega_{m}=1000\max\left\{ \left(\omega_{C}^{M}+\omega_{e}^{M}\right),\,2\omega_{C}^{M}\right\} $
implying $\frac{\omega_{m}}{2\omega_{C}^{M}}\geq1000$. For $\omega_{2}^{M}R_{f}C_{f}\leq0.1$
attenuation is more than $-20\, dB$. As the HPF offers high attenuation
to the maximum component frequency of $V_{m}^{2}$, it will offer
higher attenuation to the lower component frequency of $V_{m}^{2}$.
Hence $V_{m}^{2}$ gets highly attenuated.

As $V_{m}^{1}$ undergoes almost no attenuation ($\approx0dB$), the
output of HPF is $V_{f}=-V_{m}^{1}$. The minus sign before $V_{m}^{1}$
is justified as the HPF is in inverting mode.

\subsection{Envelope Detector\label{sub:Envelope-Detector}}

The input to this block is $V_{f}=\frac{MV_{e}\left(t\right)}{R_{I}}\sin\left(\omega_{m}t\right)$.
Similar to amplitude modulation (AM), here $\sin\left(\omega_{m}t\right)$
is the carrier and $\frac{MV_{e}\left(t\right)}{R_{I}}$ is the signal
to be recovered. We use a polarity sensitive envelope detector as
$\frac{MV_{e}\left(t\right)}{R_{I}}$ can be both positive or negative.
The key idea used here is that the polarity of $V_{e}$ and $V_{u}$
is same since $K=\frac{M}{R_{I}}>0$. Hence we detect the positive
peaks of $V_{f}$ when $V_{e}$ is positive by keeping $T_{1}$ $ON$
and $T_{C1}$ $OFF$. When $V_{e}$ is negative, negative peaks of
$V_{f}$ are detected by keeping $T_{C1}$ $ON$ and $T_{1}$ $OFF$.
Remaining working of the envelope detector is similar to a conventional
Diode-based Envelope Detector and can be found in \cite{haykin2009communication}.
Effective envelope detection using diode based envelope detector requires:

\begin{enumerate}
\item $\omega_{m}\gg\omega_{C}^{M}+\omega_{e}^{M}$, i.e. the frequency of the carrier should be much greater than the maximum component frequency of the signal. Here $\sin\left(\omega_{m}t\right)$ is the carrier whose frequency is $\omega_{m}$ while $\frac{MV_{e}}{R_{I}}$ is the signal whose maximum component frequency is $\omega_{C}^{M}+\omega_{e}^{M}$ (refer equation (8)).
\item $\frac{1}{\omega_{m}}\ll R_{e}C_{e}$, i.e. the time constant of the envelope detector is much larger than the time period of the modulating signal. This ensures that the capacitor $C_{e}$ (refer Fig. 3) discharges slowly between peaks thereby reducing the ripples.
\item $R_{e}C_{e}<\frac{1}{\omega_{C}^{M}+\omega_{e}^{M}}$, i.e. the time constant of the envelope detector should be less than the time period corresponding to the maximum component frequency of the signal getting modulated. This is necessary so that the output of the envelope detector can effectively track the envelope of the modulated signal.
\end{enumerate}

The proposed tuning rule in \textit{Claim 1} satisfies these conditions.
In general, for amplitude modulation (AM) the choice of modulating
frequency is $100$ times the maximum component frequency of the signal
getting modulated. Unlike AM, our multiplying factor is $1000$ instead
of $100$. The reason for choosing this should be clearly understood.
In a conventional diode based peak detector the ripple in output voltage
can be decreased either by increasing the modulating frequency $\omega_{m}$
or by increasing the time constant $R_{e}C_{e}$. But $R_{e}C_{e}$
is upper bounded by $\frac{1}{\omega_{C}^{M}+\omega_{e}^{M}}$. In
our case the signal to be modulated, i.e. $\frac{MV_{e}}{R_{I}}$,
may contain frequencies anywhere in the range $\left[0,\omega_{C}^{M}+\omega_{e}^{M}\right]$.
For a given $R_{e}C_{e}$ a signal with a lower frequency will suffer
higher ripple. Hence to constrain the ripple for any frequency in
the specified range one must constrain the ripple for $0$ frequency
(DC voltage). For DC voltage ripple factor is approximately $\frac{2\pi\times100}{\sqrt{3}\omega_{m}R_{e}C_{e}}$.
With the choice of $R_{e}C_{e}$ made in Claim 1 and $100$ as multipying
factor, ripple factor is as large as $7.2\%$. Therefore, we choose
multiplying factor of $1000$ which gives a ripple factor of $0.72\%$.
The output of the envelope detector is $\frac{MV_{e}\left(t\right)}{R_{I}}$
and the final output of the circuit is

\vspace{-0.5em}

\begin{equation}
\dot{\mathcal{M}}\thickapprox-\frac{\alpha^{S}}{R_{I}R_{C}}V_{C}\left(t\right)\,;\; V_{u}\left(t\right)=\mathcal{M}V_{e}
\end{equation}
where, $\mathcal{M}=\frac{M}{R_{I}}$. Comparing equations (5) and
(12) we see that $\alpha_{k}=-\frac{\alpha^{S}}{R_{I}R_{C}}$ and
$K=\mathcal{M}$ where $\mathcal{M}\in\left[\frac{R_{on}^{S}}{R_{I}}\,,\:\frac{R_{off}^{S}}{R_{I}}\right]$.
$R_{I}$ is tuned to get the desired range of gain while $R_{C}$
is a free parameter which can be tuned according to the needs.

\subsection{Charge Saturator Block\label{sub:Charge-Saturator-Block}}

This block limits the memristor to work in its safe zone hence ensuring
validity of equation (4). In safe zone the following equations are
valid%
\footnote{Actually, $\frac{dQ_{M}}{dt}=\frac{V_{C}^{m}}{R_{C}}+\frac{V_{e}\sin\left(\left(\omega_{m}t\right)\right)}{R_{I}}$.
But as proved in \ref{sub:Memristor-Gain-Block}, the effect of high
frequency term, $V_{e}\sin\left(\omega_{m}t\right)$ on $Q_{M}$ is
negligible.%
},

\vspace{-0.8em}

\begin{equation}
\frac{dV_{IG}}{dt}=-\frac{V_{C}^{m}}{R_{s}C_{s}}\:;\:\frac{dQ_{M}}{dt}=\frac{V_{C}^{m}}{R_{C}}\:\Rightarrow\:\frac{dV_{IG}}{dQ_{M}}=-\frac{R_{C}}{R_{s}C_{s}}
\end{equation}

Recall that in the safe zone $w\in\left[w_{l}\,,\, w_{h}\right]$
and $Q_{M}\in\left[0\,,\, Q_{M}^{S}\right]$. We assume that integrator
voltage $V_{IG}=0$ when $Q_{M}=0$ (or $w=w_{l}$). Integrating equation
(13) under this assumption yields

\vspace{-0.6em}

\begin{equation}
V_{IG}=-\frac{R_{C}}{R_{s}C_{s}}Q_{M}
\end{equation}

\vspace{-0.6em}

\noindent In equation (14), if $Q_{M}=Q_{M}^{S}$ (or $w=w_{h}$),
$V_{IG}=-\frac{R_{C}Q_{M}^{S}}{R_{s}C_{s}}$. Hence, $V_{IG}\in\left[-\frac{R_{C}Q_{M}^{S}}{R_{s}C_{s}}\,,\,0\right]$
when memristor is in its safe zone. In Fig. 3 comparator $O_{1}$
and $O_{2}$ are used to compare $V_{IG}$ to know if the memristor
is in safe zone. Note that: 1) Reference voltage of comparator $O_{1}$
and $O_{2}$ is $GND$ and voltage $V_{H}$ (across capacitor $C_{H}$)
respectively. $V_{H}$ is set to $-\frac{R_{C}Q_{M}^{S}}{R_{s}C_{s}}$
by Synchronization Block (refer Section \ref{sub:SynchroBlock}).
2) In Fig. 3 any MOSFET transistor couple, $T_{i}$ and $T_{Ci}$,
are in complementary state, i.e. if $T_{i}$ is $ON$, $T_{Ci}$ will
be $OFF$ and vice-versa. 3) Comparator output $V_{L1}$ and $V_{L2}$
gives knowledge about the state of the memristor. Also $V_{L1},V_{L2}\in\left\{ -V_{DD},V_{DD}\right\} $.
Now depending on $V_{L1}$ and $V_{L2}$, three different cases may
arise:

\begin{figure}[t]
\begin{centering}
\includegraphics[width=8.9cm,height=7cm]{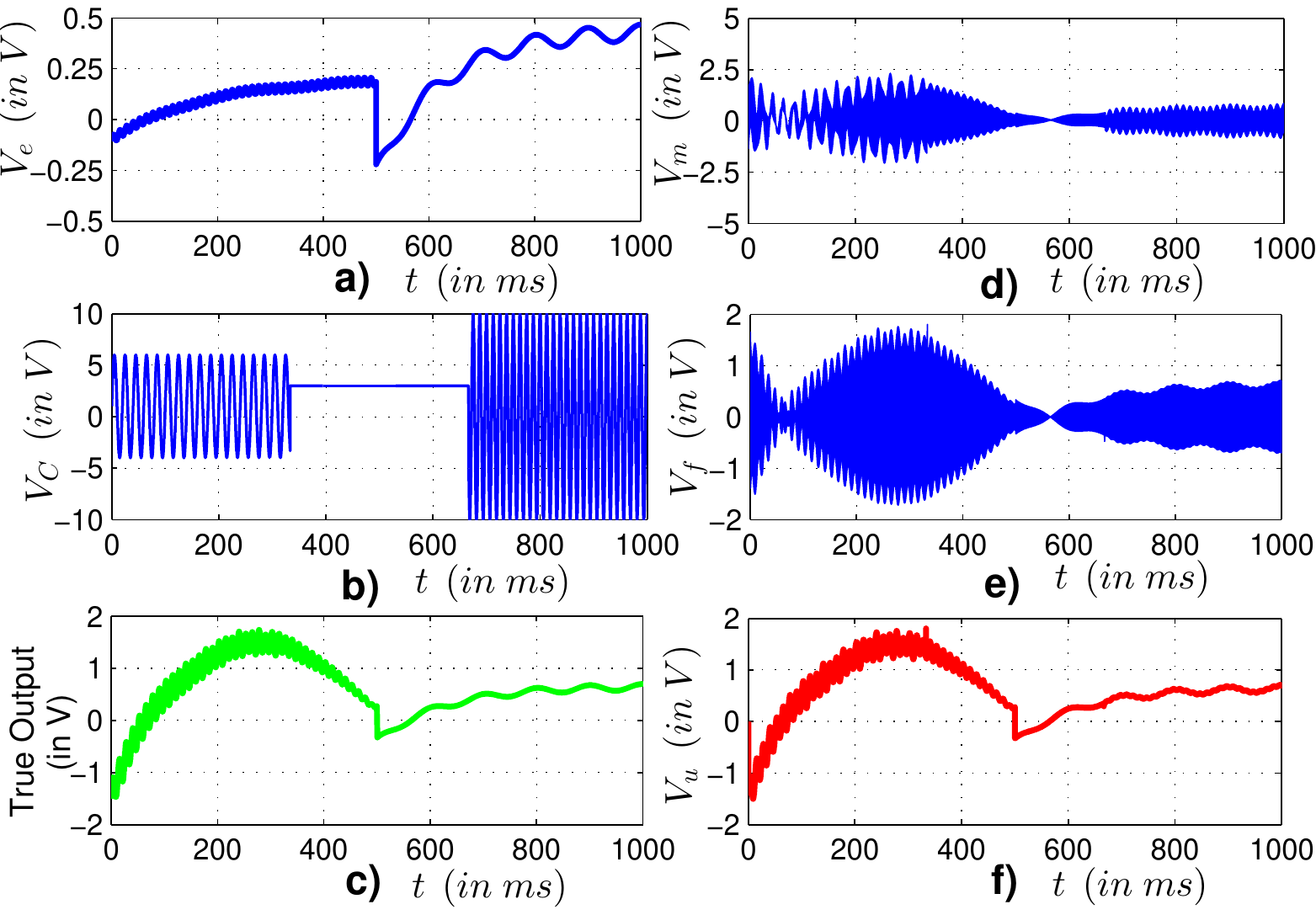} \caption{Output of various stages of the circuit shown in Fig. 3, i.e. $V_{m}$,
$V_{f}$ and $V_{u}$, corresponding to inputs: $V_{e}$ and $V_{C}$.
Parameters of simulations are: HP Memristor in safe zone, $\omega_{m}=628\times10^{3}$
$rad s^{-1}$, $R_{C}=100$ $k\Omega$, $R_{I}=1$ $k\Omega$, $R_{f}C_{f}=0.159$
$ms$, $R_{e}C_{e}=0.796$ $ms$, $V_{DD}=5$ $V$, $R_{s}C_{s}=0.826$
$s$. $V_{u}$ is obtained by simulating the circuit show in Fig.
3. The $true$ $output$ is obtained by numerically solving equation
(5). In both cases we use $V_{e}$ and $V_{C}$ shown in Fig. 4a and
Fig. 4b respectively as inputs.}

\par\end{centering}

\vspace{-1.0em}
\end{figure}

\textbf{Case 1 }($V_{H}<V_{IG}<0\Rightarrow\ensuremath{V_{L1}=V_{DD}}\,,\: V_{L2}=V_{DD}$)

\noindent This implies that the memristor is in its safe zone, i.e.
$0<Q_{M}<Q_{M}^{S}$. $Q_{M}$ can either increase or decrease. Hence
both $T_{2}$ and $T_{3}$ are $ON$ allowing both positive and negative
voltage to pass through, i.e. $V_{C}^{m}=V_{C}$.

\textbf{Case 2 }($V_{IG}<V_{H}\Rightarrow\ensuremath{V_{L1}=V_{DD}}\,,\: V_{L2}=-V_{DD}$)

\noindent This happens when $Q_{M}\geq Q_{M}^{S}$. Since $Q_{M}$
can only decrease, $T_{2}$ is kept $OFF$ but $T_{3}$ is $ON$.
Two cases are possible: if $V_{C}>0$, $T_{4}$ and $T_{C2}$ will
be $ON$ making $V_{C}^{m}=0$ or $V_{C}<0$ making it pass through
$T_{3}$ and thereby setting $V_{C}^{m}=V_{C}$.

\textbf{Case 3 }($V_{IG}>0\Rightarrow\ensuremath{V_{L1}=-V_{DD}}\,,\: V_{L2}=V_{DD}$)

\noindent This happens when $Q_{M}\leq0$. Since $Q_{M}$ can only
increase, $T_{3}$ is kept $OFF$ but $T_{2}$ is $ON$. Two cases
are possible: if $V_{C}<0$, $T_{C4}$ and $T_{C3}$ will be $ON$
making $V_{C}^{m}=0$ or $V_{C}>0$ making it pass through $T_{2}$
and thereby setting $V_{C}^{m}=V_{C}$.

\vspace{-0.3em}

\subsection{Synchronization Block\label{sub:SynchroBlock}}

Operation of Charge Saturator Block assumes that: 1) $V_{IG}=0$ when
$Q_{M}=0$ (or $w=w_{l}$). 2) Voltage $V_{H}$ across capacitor $C_{H}$
equals $-\frac{R_{C}Q_{M}^{S}}{R_{s}C_{s}}$. This block ensures these
two conditions and thereby guarantees that the memristor and the integrator
in Fig. 3 are in ``synchronization''. Synchronization Block is shown
in Fig. 5a). In Fig. 5a, the op-amp with the memristor, the integrator
and capacitor $C_{H}$ are indeed part of the circuit shown in Fig.
3. Such a change in circuit connection is possible using appropriate
switching circuitry. This block operates in two modes:

\begin{figure}[t]
\begin{centering}
\includegraphics[scale=0.36]{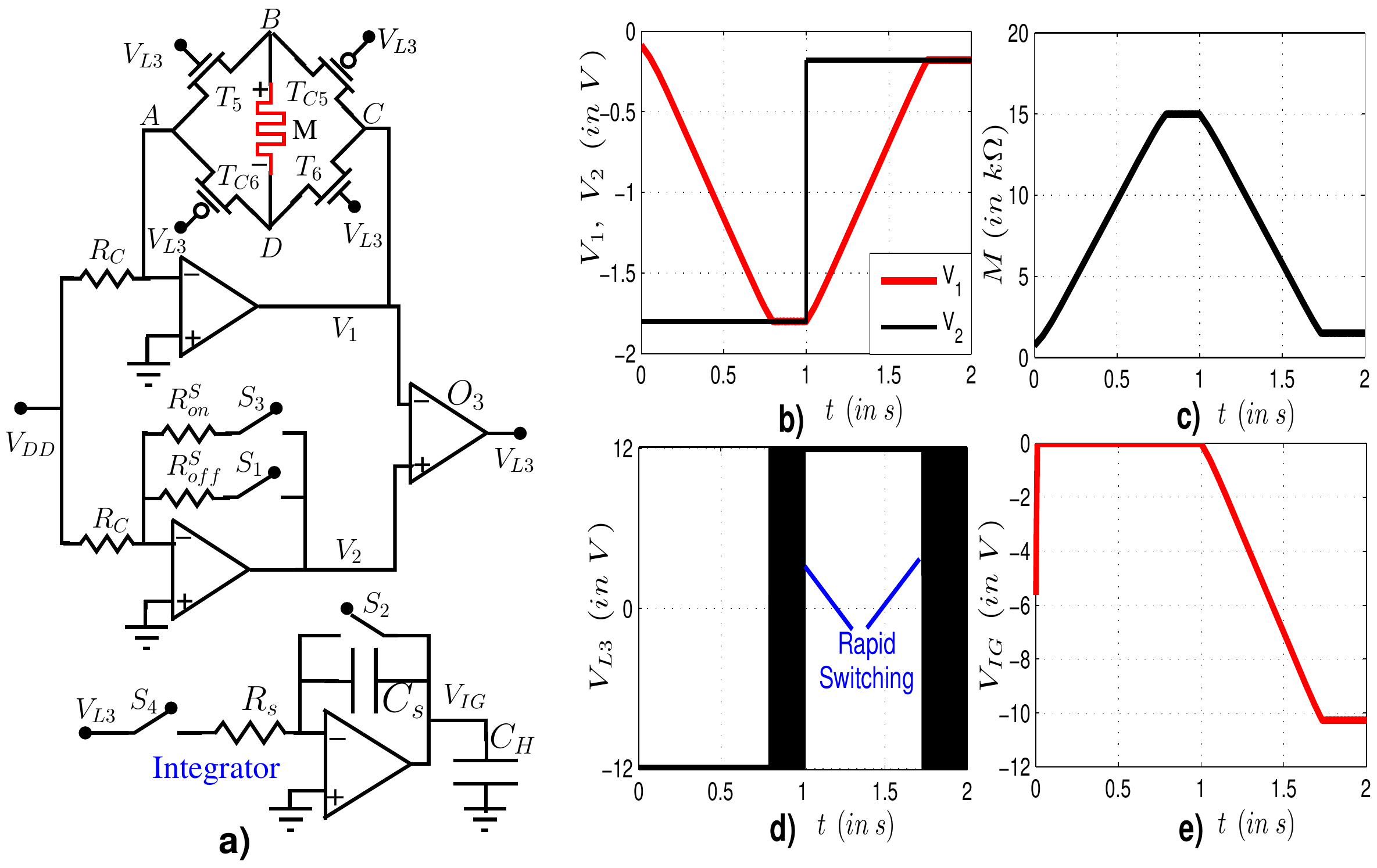}
\par\end{centering}

\caption{a) Schematic of Synchronization Block. b), c), d), e) Graphs showing operations of Synchronization Block. In these graphs preset mode operates for the first $1s$ and online calibration mode operates for the next $1s$.}

\vspace{-0.8em}
\end{figure}

\noindent \textbf{Preset Mode: }We first ensure that $V_{IG}=0$,
when $w=w_{l}$. In this mode switch $S_{1}$ and $S_{2}$ are $ON$
and switch $S_{3}$ and $S_{4}$ are $OFF$. When $S_{2}$ is closed
the residual charge in capacitor $C_{H}$ will get neutralized thus
ensuring $V_{IG}=0$. Next we make $w=w_{l}$. Note that $V_{1}=-\frac{MV_{DD}}{R_{C}}$
and $V_{2}=-\frac{R_{off}^{S}V_{DD}}{R_{C}}$. If $M<R_{off}^{S}$
then $V_{1}>V_{2}\,\Rightarrow V_{L3}=-V_{DD}$. Hence the path ADBC
of wheatstone bridge arrangement will be active making the current
flow from $\left(-\right)^{ve}$ to $\left(+\right)^{ve}$ terminal
of the memristor. This will increase $M$ till $M=R_{off}^{S}$. If
$M>R_{off}^{S},$ the path ABDC is active making $M$ decrease till
$M=R_{off}^{S}$.

\noindent \textbf{Online Calibration\footnote{If we tune $R_{s}$ and $C_{s}$ s.t. $-\frac{R_{C}Q_{M}^{S}}{R_{s}C_{s}}=-V_{DD}$, $V_{IG}\in\left[-V_{DD},0\right]$ when memristor is in safe zone. Then we can directly use power supply $-V_{DD}$ as the reference voltage for comparator $O_{2}$. This will eliminate the need of capacitor $C_{H}$ and ``Online Calibration''. However if $-\frac{R_{C}Q_{M}^{S}}{R_{s}C_{s}}\neq-V_{DD}$ (due to tuning error), this approach may drive the memristor to non-safe zone.}: }Immediately
after preset mode is complete, $S_{1}$ and $S_{2}$ are switched
$OFF$ and $S_{3}$ and $S_{4}$ are switched $ON$. Now, $V_{1}=-\frac{MV_{DD}}{R_{C}}$
and $V_{2}=-\frac{R_{on}^{S}V_{DD}}{R_{C}}$. In this step $V_{L3}=V_{DD}$
as $M>R_{on}^{S}$ always. Path ABDC will be active driving $M$ to
$R_{on}^{S}$. As $V_{L3}$ is given as an input to the integrator,
capacitor $C_{H}$ will also get charged. Note that in this step memristor
will work in safe zone. Also $V_{IG}=0$ when $Q_{M}=0$ (ensured
by preset mode). Hence relation between $V_{IG}$ and $Q_{M}$ will
be governed by equation (14). When $M$ gets equalized to $R_{on}^{S}$,
$Q_{M}=Q_{M}^{S}$, thereby making $V_{H}=V_{IG}=-\frac{R_{C}Q_{M}^{S}}{R_{s}C_{s}}$.

Each of the modes operate for a predefined time. The resistors and
hence the voltages $V_{1}$ and $V_{2}$ may get equalized before
the predefined time after which $V_{L3}$ will switch rapidly. Such
rapid switching can be prevented by replacing the comparator $O_{3}$
by a cascaded arrangement of a differential amplifier followed by
a hysteresis block. 

Various graphs corresponding to synchronization process are shown
in Fig. 5 b), c), d), e). Memristor Gain Block and the Integrator
of Charge Saturator Block should be periodically synchronized to account
for circuit non-idealities. One such non-ideality can be caused if
the capacitor $C_{H}$ is leaky causing the voltage $V_{H}$ to drift
with time. 

\textit{Remark 3: }In the discussion of the Synchronization Block
we have slightly misused the symbols $R_{on}^{S}$ and $R_{off}^{S}$.
Resistance of $R_{on}^{S}$ and $R_{off}^{S}$ shown in Fig. 5a should
be close to the actual $R_{on}^{S}$ and $R_{off}^{S}$ (as mentioned
in Section \ref{sec:Prelim}) respectively. It is not necessary that
they should be exactly equal. However, there resistance must lie within
the safe zone. This alleviates analog implementation by eliminating
the need of \textit{precision resistors}. It should be noted that
the maximum and the minimum resistance of the memristor in safe zone
is governed by the resistances $R_{on}^{S}$ and $R_{off}^{S}$ used
in the Synchronization Block \textit{not} the actual $R_{on}^{S}$
and $R_{off}^{S}$.

To conclude, in this section we designed an Analog Gain Control framework
using Memristor. Schematic of Memristive AGC is shown in Fig. 2 whose
circuit parameters can be tuned using \textit{Claim 1}. Memristive
AGC's designed in this work is ``generic'' in the sense that it
can be used to implement several Gain-Scheduled control algorithms.

\section{Control Strategy\label{sec:cntrl}}

As mentioned in the introduction, we are interested in designing Robust
Adaptive Control Algorithms to tackle issues like parameter uncertainty
(caused by design errors) and time varying nature (caused by ageing
effect and atmospheric variation). `\textit{Simplicity}' is the
key aspect of any control algorithm to be implementable in analog
framework as we do not have the flexibility of `\textit{coding}'.
Robust Adaptive Control Algorithms found in control literature cannot
be implemented in an analog framework due to their complexity. Here
we propose a simple gain-scheduled robust adaptive control algorithm
which can be easily implemented using Memristive AGC discussed in
Section \ref{sec:memgain}. We prove the stability of the proposed
algorithm using \textit{Lyapunov-Like} method in Section \ref{sub:Reflective-Gain-Space}.

\textit{Notations: }The notations used in this paper are quite standard.
$\mathbb{R}^{+}$ and $\mathbb{R}^{n}$ denotes the set of positive
real numbers and the n-dimensional real space respectively. $I$ represents
identity matrix of appropriate dimension. $\left|\cdot\right|$ is
the absolute value operator. $\emptyset$ represents a null set. The
bold face symbols $\boldsymbol{S}$ and $\mathbf{S}^{+}$ represents
the set of all symmetric matrices and positive definite symmetric
matrices respectively. $\inf\left(\cdot\right)$ ($\sup\left(\cdot\right)$)
represents the infimum (supremum) of a set. For a matrix $A$, $\lambda_{m}\left(A\right)$
and $\lambda_{M}\left(A\right)$ denotes the minimum and maximum eigenvalue
of $A$ respectively. $A\preceq B$ implies $B-A$ is positive semi-definite
while $A\prec B$ implies $B-A$ is positive definite. The euclidean
norm of a vector and the induced spectral norm of a matrix is denoted
$\begin{Vmatrix}\cdot\end{Vmatrix}$. The operator $\times$ when
applied on sets implies the cartesian product of the sets. $Conv\left\{ \mathcal{A}\right\} $
implies the convex hull of set $\mathcal{A}$. Analysis presented
futher in this paper uses ideas from real analysis, linear algebra
and convex analysis. For completeness, these ideas are briefly reviewed
in Appendix \ref{sec:Notions}.

\begin{figure}[t]
\begin{centering}
\includegraphics[scale=0.3]{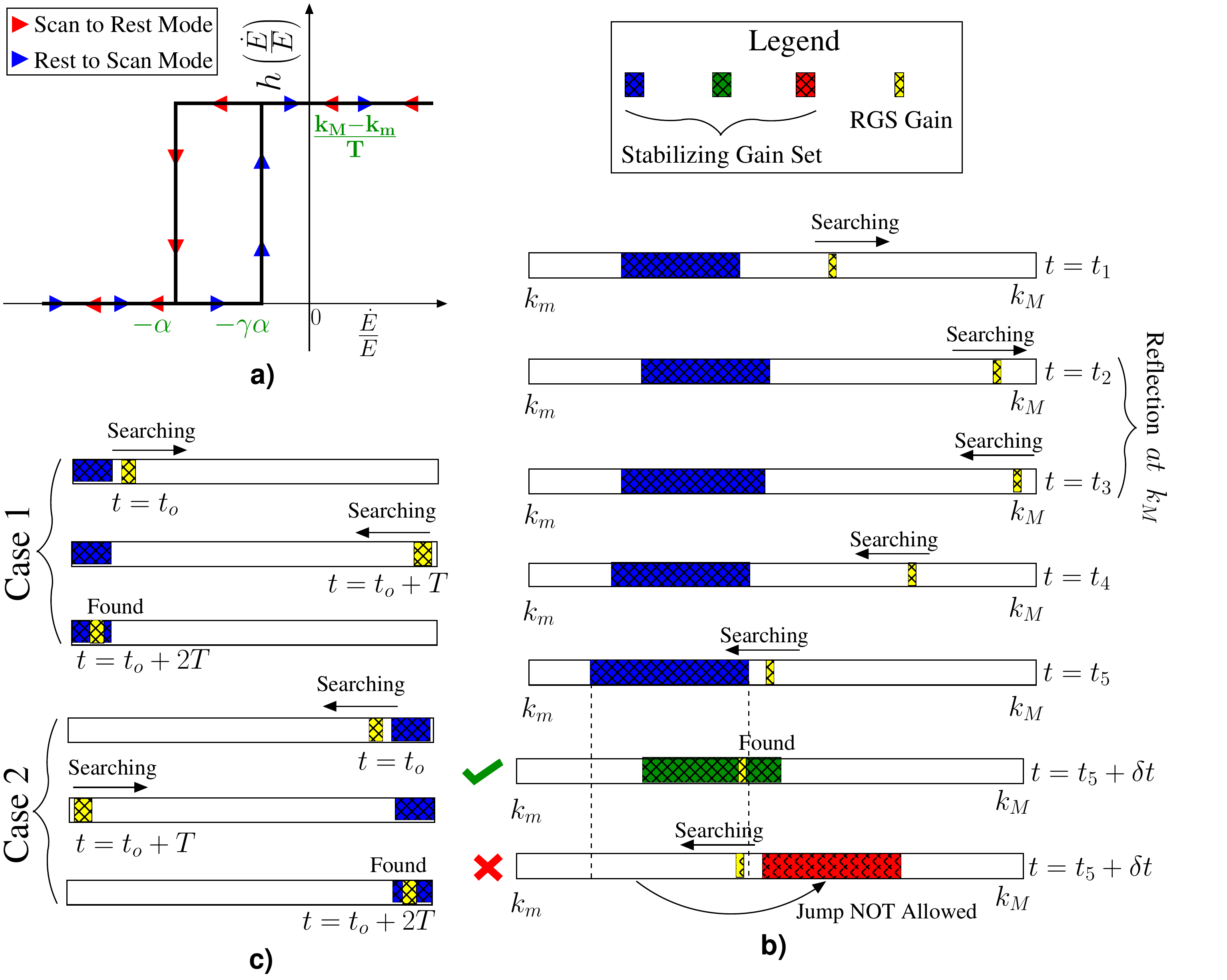}
\par\end{centering}

\caption{a) Hysteresis Function as mentioned in equation (16). $T$ is the scan time. $\alpha\in\mathbb{R^{+}}-\left\{ 0\right\}$ and $\gamma\in\left[0,1\right)$ are tuning constants. b) First six images shows the scan/rest mode of RGS. The last three images depicts the concept of \textit{drifting}. c) Figure showing the worst case scan time.}

\vspace{-1.0em}
\end{figure}

\subsection{Reflective Gain Space Search (RGS)\label{sub:Reflective-Gain-Space}}

Consider a SISO uncertain linear time variant (LTV) system with states
$x=\begin{bmatrix}x_{1} & x_{2} & \cdots & x_{N}\end{bmatrix}^{T}\in\mathbb{R}^{N}$,
input $u\in\mathbb{R}$ and output $y\in\mathbb{R}$ described by
the following state space equation

\vspace{-1.2em}

\begin{equation}
\dot{x}=A\left(t\right)x+B\left(t\right)u\;;\quad y=Cx=x_{1}\quad\text{{where,}}
\end{equation}

$A\left(t\right)=\left[\begin{array}{cccc}
0 & 1 & \cdots & 0\\
\vdots & \vdots & \ddots & \vdots\\
0 & 0 & \cdots & 1\\
a_{1}\left(t\right) & a_{2}\left(t\right) & \cdots & a_{N}\left(t\right)
\end{array}\right]\; B\left(t\right)=\begin{bmatrix}0\\
\vdots\\
0\\
b\left(t\right)
\end{bmatrix}$

The output matrix $C=\left[\begin{array}{cccc}
1 & \cdots & 0 & 0\end{array}\right]$.

\textit{Assumptions}:

\begin{enumerate}
\item The order $N$ of the system is known. 
\item The variation of $a_{i}\left(t\right),\,\forall i=1,2,\ldots,N$ and $b\left(t\right)$ due to parameter uncertainty and time variation is bounded, i.e. $\left(A\left(t\right),B\left(t\right)\right)$ belongs to a bounded set $\mathcal{L}$. $\mathcal{L}$ is assumed to be a connected set. $A\left(t\right)$ and $B\left(t\right)$ are not known but knowledge of $\mathcal{L}$ is assumed. This assumption basically means that we don't have knowledge about $a_{i}\left(t\right)$ and $b\left(t\right)$ however we know the \textit{range} in which they lie.
\item $\begin{Vmatrix}\dot{A}\end{Vmatrix}\leq\delta_{A}$ and $\begin{Vmatrix}\dot{B}\end{Vmatrix}\leq\delta_{B}$. Knowledge of $\delta_{A}$ and $\delta_{B}$ is assumed. This assumption basically means that $A\left(t\right)$ and $B\left(t\right)$ are continous in time.
\item Let $b\left(t\right)\in\left[b_{l},b_{u}\right]$ s.t. $b_{l}\leq b_{u}$, and either $b_{l}>0$ or $b_{u}<0$. This choice of $b_{l}$ and $b_{u}$ is explained in Section~\ref{sub:optimization}.
\end{enumerate}

Example to clarify the concept of $\mathcal{L}$, $\delta A$ and
$\delta B$ is dealt later in Section \ref{sub:optimization} however
we give the following example to better explain \textit{Assumption 4}.

\textit{Example 1: }Consider two cases: 1) $b\left(t\right)\in\left[-2,\,3\right]$
2) $b\left(t\right)\in\left[0.5,\,3\right]$. According to \textit{Assumption 4},
\textit{Case 2} is possible while \textit{Case 1} is not. This example
clearly illustrates that according to \textit{Assumption 4} the sign
of $b\left(t\right)$ does not change with time and is known with
absolute certainty. Note that sign of $b\left(t\right)$ decides the
sign of static gain%
\footnote{Unlike LTI system, static gain for LTV system may not be well defined.
Here we define static gain of LTV system (15) as $\frac{b\left(t\right)}{a_{1}\left(t\right)}$.%
} of the system. So \textit{Assumption 4} in certain sense means that
the sign of static gain does not change with time and is known with
absolute certainty. For all practical scenario such an assumption
is valid.

Notice that all assumptions are mild from practical viewpoint. Our
aim is to regulate the output around the operating point $y^{*}=x_{1}^{*}=r$.
Conventional set-point control consist of a bias term $u_{b}\left(t\right)$
plus the regulation control input $u_{r}\left(t\right)$, i.e. $u\left(t\right)=u_{b}\left(t\right)+u_{r}\left(t\right)$.
Here we assume that $u_{b}\left(t\right)$ is designed s.t. $x^{*}=\begin{bmatrix}r & 0 & \ldots & 0\end{bmatrix}^{T}$
is the equilibrium point and concentrate on the synthesis of $u_{r}\left(t\right)$.
For simplicity we consider $r=0$, i.e. $x^{*}=\begin{bmatrix}0 & 0 & \ldots & 0\end{bmatrix}^{T}$
is the equilibrium point. As we are dealing with a linear system the
same analysis is valid for $r\neq0$. The controller structure is,
$u_{r}\left(t\right)=-K\left(t\right)y\left(t\right)$, where $K\left(t\right)\in\left[k_{m},k_{M}\right]$
is the variable gain s.t. $0<k_{m}<k_{M}$. Let $E=x^{T}Px$ be the
Lyapunov Candidate Function%
\footnote{Use of time invariant Lyapunov Function to analyse stability of LTV
systems has been used in \cite{boyd1994linear} (Chap. 5, 7) and \cite{ben2001lectures}
(Chap. 3).%
} with $P\in\mathbf{S}^{+}$. Then RGS is as simple Gain-Scheduled
control strategy given by the following equation

\vspace{-0.6em}

\begin{equation}
\dot{K}=sgn\cdot h\left(\frac{\dot{E}}{E}\right)\;;\quad u_{r}=-Ky
\end{equation}

\vspace{-0.6em}

\noindent where, $sgn\in\left\{ -1,1\right\} $ and $h\left(\frac{\dot{E}}{E}\right)$
is a hysteresis function shown in Fig. 6a. Working of RGS can be explained
as:

\noindent \begin{enumerate}
\item RGS finds the stabilizing gain\footnote{The term ''stabilizing gain'' has been slightly misused. Stability of LTV system cannot be assured even if the closed loop system matrix, $A\left(t\right)-B\left(t\right)K\left(t\right)C$, has negative real eigen part for all $t>0$ (\cite{rosenbrook,solo}).}, i.e. the gain which renders $\dot{E}<0$ ($\dot{E}<-\alpha E$ in a more strict sense), by reflecting back and forth between $\left[k_{m},\: k_{M}\right]$. RGS is said to be in "Scan Mode" when it is scanning for the stabilizing gain. It goes to "Rest Mode" when stabilizing gain is found. RGS Scan Cycle is clearly depicted in the first six images of Fig. 5b).
\item RGS uses $\dot{E}$ as stability indicator. $\dot{E}$ is found by differentiating $E$ which in turn is calculated using $E=x\left(t\right)^{T}Px\left(t\right)$. To get the states $x\left(t\right)$ we differentiate the output $y\left(t\right)$ $N-1$ times.
\item Scan Mode is triggered when  $\dot{E} > - \gamma \alpha E$ (refer Fig. 6a). Scan Mode is associated with a scan time of $T$, i.e. time taken to scan from $k_{m}$ to $k_{M}$. Hence, $h\left(\frac{\dot{E}}{E}\right)=\frac{k_{M}-k_{m}}{T}$. The value of $sgn$ is $1$ when gain space is searched from $k_{m}$ to $k_{M}$ and $-1$ otherwise. Scan mode operates till $\dot{E}>-\alpha E$.
\item Rest Mode is triggered when $\dot{E} < - \alpha E$. In this mode $h\left(\frac{\dot{E}}{E}\right)=0$, i.e. the stabilizing gain is held constant. Rest mode operates till $\dot{E} < - \gamma \alpha E$.
\item In the process of finding stabilizing gain, LTV system \textit{may} expand ($E$ increases) in Scan Mode and \textit{will} contract ($E$ decreases) in Rest Mode. RGS ensures that even in the worst case, contraction is always dominant over expansion, guaranteeing stability.
\end{enumerate}

Stabilizing Gain Set%
\footnote{The stabilizing gain set $\mathcal{K}_{P,s}\left(t_{o}\right)$ at
time $t=t_{o}$, is just a ``Lyapunov-Way'' of describing the set
of gains which will stabilize the corresponding LTI system $\left(A\left(t_{o}\right),B\left(t_{o}\right),C\right)$
at time $t=t_{o}$.%
} $\mathcal{K}_{P,s}\left(t\right)$ for a given $P\succ0$ and $s\in\mathbb{R}^{+}$
at time $t$ is defined as

\vspace{0.5em}

$\mathcal{K}_{P,s}\left(t\right)=\left\{ K\in\mathbb{R}\,:\: A_{C}\left(t\right)^{T}P+PA_{C}\left(t\right)\preceq-sI,\right.$

$\qquad\qquad\quad\left.\vphantom{A_{C}\left(t\right)^{T}}A_{C}\left(t\right)=A\left(t\right)-B\left(t\right)KC,\: k_{m}\leq K\leq k_{M}\right\} $

\vspace{0.5em}

\begin{lemma}
If $\mathcal{K}_{P,s}\left(t\right)\neq\emptyset,\;\forall t\geq0$ then under \textit{Assumption 3}, $\mathcal{K}_{P,s}\left(t\right)\bigcap\mathcal{K}_{P,s}\left(t+\delta t\right)\neq\emptyset$ if $\delta t \rightarrow 0$.
\end{lemma}

\textit{Proof: }This lemma basically means that stabilizing gain
set will drift. If the stabilizing gain set is drifting then there
will be an intersection between stabilizing gain sets at time $t_{1}$
and $t_{2}$ if $t_{2}-t_{1}$ is small. Drifting is obvious under
\textit{Assumption 3} however a formal proof is given in Appendix
\ref{sec:Proof-Lemma1}.

Concept of ``drifting'' is clearly depicted in the last three images
of Fig. 5b. Notice that at $t=t_{5}$ the stabilizing gain set is
almost found. Two possible cases are shown in Fig. 5b for $t=t_{5}+\delta t$
where $\delta t\rightarrow0$. In the first case the stabilizing gain
set is drifting while in the other case the stabilizing gain set jumps.
In the first case the stabilizing gain is found and RGS goes to Rest
Mode. In the second case RGS will just miss the stabilizing gain set.
Hence if the stabilizing set keeps jumping, the scan mode may never
end. Hence if the stabilizing set keeps jumping, the scan mode may
never end. As mentioned in \textit{Lemma 1}, stabilizing set $\mathcal{K}_{P,s}\left(t\right)$
never jumps if \textit{Assumption 3} is valid. Therefore \textit{Lemma 1}
is important to guarantee an upper bound on the time period of Scan
Mode from which we get the following Lemma.

\begin{lemma}
If $\mathcal{K}_{P,s}\left(t\right)\neq\emptyset,\;\forall t\geq0$ then under \textit{Assumption 3}, the maximum time period of Scan Mode is $2T$.
\end{lemma}

\textit{Proof: }It is a direct consequence of \textit{Lemma 1}.
The worst case time period of $2T$ happens only in two cases. Both
the cases are shown in Fig. 6c.

\begin{theorem}
LTV system characterized by equation (15) and set $\mathcal{L}$ is stable under RGS Control Law proposed in equation (16) if it satisfies the following criterion
\begin{enumerate}
\item \textbf{[C1] }If corresponding to all $\left(A,B\right)\in\mathcal{L}$ there exist atleast one gain $K_{AB}\in\left[k_{m},k_{M}\right]$ s.t.
\begin{equation} \left(A-BK_{AB}C\right)^{T}P+P\left(A-BK_{AB}C\right)\preceq-sI \end{equation}
for some $s\in\mathbb{R^{+}}$.
\item \textbf{[C2] } Scan time $T<\frac{\lambda_{m}\left(P\right)\alpha^{2}\left(1-\gamma^{2}\right)}{8\delta\max\left(\lambda_{M}^{\mathcal{L}},0\right)}$. Here $\gamma\in\left[0,1\right)$, $\delta=\delta_{A}+\delta_{B}k_{M}$ and $\lambda_{M}^{\mathcal{L}}=\underset{Q\in S_{Q}^{\mathcal{L}}}{\max}\lambda_{M}\left(Q\right)$ where the set
\end{enumerate}
\end{theorem}

$\quad S_{Q}^{\mathcal{L}}=\left\{ Q\,:\: Q=\left(A-BKC\right)^{T}P+P\left(A-BKC\right)\right.$

$\quad\quad\quad\;\quad\left.\vphantom{A_{C}\left(t\right)^{T}}\left(A,B\right)\in\mathcal{L},\: K\in\left[k_{m},k_{M}\right]\right\} $

\textit{Proof: }Consider the Lyapunov candidate $E=x^{T}Px$ which
when differentiated along systems trajectory yields

\vspace{-1.2em}

\begin{equation}
\dot{E}=x^{T}Q\left(t\right)x\leq\lambda_{M}\left(Q\left(t\right)\right)\begin{Vmatrix}x\end{Vmatrix}^{2}\leq\lambda_{M}^{\mathcal{L}}\begin{Vmatrix}x\end{Vmatrix}^{2}
\end{equation}

\vspace{-0.3em}

\noindent where, $Q\left(t\right)=A_{C}\left(t\right)^{T}P+PA_{C}\left(t\right)$
and $A_{C}\left(t\right)=A\left(t\right)-B\left(t\right)K\left(t\right)C$.
Definitely, $Q\left(t\right)\in\mathbf{S}$ as it is a sum of the
matrix $PA_{C}\left(t\right)$ and its transpose $A_{C}\left(t\right)^{T}P$.
Note that, $\lambda_{M}\left(Q\left(t\right)\right)\leq\lambda_{M}^{\mathcal{L}},\:\forall\left(A,B\right)\in\mathcal{L}\:\text{{and}}\:\forall K\in\left[k_{m},k_{M}\right]$.
When $\lambda_{M}^{\mathcal{L}}<0$, stability is obvious as according
to inequality (18) $\dot{E}<0$, $\forall\left(A,\, B\right)\in\mathcal{L}$
and $\forall K\in\left[k_{m},\, k_{M}\right]$. As $\dot{E}<0$ for
any $\forall K\in\left[k_{m},k_{M}\right]$, stability if guaranteed
even if scanning is infinitely slow, i.e. scan time $T\rightarrow\infty$.
This is in accordance with \textbf{[C2]}.

We now consider the case when $\lambda_{M}^{\mathcal{L}}>0$. Say
at time $t=t^{*}$, $E=E^{*}$ and the system goes to scan mode to
search for a stabilizing gain. From inequality (18) we have,

\vspace{-1.2em}

\begin{eqnarray}
\dot{E} & \leq & \lambda_{M}^{\mathcal{L}}\begin{Vmatrix}x\end{Vmatrix}^{2}\leq\frac{\lambda_{M}^{\mathcal{L}}}{\lambda_{m}\left(P\right)}E
\end{eqnarray}

\noindent We want to find the maximum possible increase in $E$ in
the Scan Mode. At this point it is important to understand that \textbf{[C1]}
is another way of saying that $\mathcal{K}_{P,s}\left(t\right)\neq\emptyset,\:\forall t\geq0$.
Hence if \textbf{[C1]} is true, then \textit{Lemma 2} can be used
to guarantee that the maximum period of scan mode is $2T$. Let $E=E_{s}$
and $t=t_{s}$ at the end of scan mode. Maximum expansion in $E$
is obtained by integrating inequality (19) from $t=t^{*}$ to $t=t^{*}+2T$ 

\vspace{-1.6em}

\begin{eqnarray}
E_{S} & \leq & \beta_{s}E^{*}\quad\text{{where,}}
\end{eqnarray}

\noindent $\beta_{s}=\exp\left(\frac{2\lambda_{M}^{\mathcal{L}}T}{\lambda_{m}\left(P\right)}\right)$,
the worst case expansion factor. The end of scan mode implies that
a stabilizing gain is found, i.e. at $t=t_{s}$, $\dot{E}\left(t_{s}\right)\leq-\alpha E\left(t_{s}\right)$
(refer Fig. 6a), and the rest mode starts. Here it is important to
note the relation between $\alpha$ and $s$ (refer \textbf{[C1]}).
\textbf{[C1]} assures that at $t=t_{s}$ a gain $K$ can be found
s.t.

\vspace{-1.5em}

\begin{eqnarray}
\dot{E}\left(t_{s}\right)\leq\lambda_{M}\left(Q\left(t_{s}\right)\right)\left\Vert x\right\Vert ^{2} & = & -s\left\Vert x\right\Vert ^{2}\nonumber \\
 & \leq & -\frac{s}{\lambda_{M}\left(P\right)}E\left(t_{s}\right)
\end{eqnarray}

\noindent Comparing inequality (21) with Fig. 6a we get $\alpha=\frac{s}{\lambda_{M}\left(P\right)}$.

Let $\tau$ be the duration of rest mode. We want to find the minimum
possible decrease in $E$ in rest mode before it again goes back to
scan mode. Again,

\vspace{-1.5em}
\begin{eqnarray}
\dot{E} & = & x^{T}Q\left(t\right)x=x^{T}\left[Q\left(t_{s}\right)+\Delta\left(t\right)\right]x\nonumber \\
 & = & \dot{E}\left(t_{s}\right)+x^{T}\Delta\left(t\right)x\quad\text{{where,}}
\end{eqnarray}

\noindent $\Delta\left(t\right)=Q\left(t\right)-Q\left(t_{s}\right)$.
As $K\left(t\right)=K\left(t_{s}\right);\, t_{s}\leq t\leq t_{s}+\tau$
($\dot{K}=0$ in rest mode), $\Delta\left(t\right)$ can be expanded
as

\vspace{-1.5em}

\begin{eqnarray}
\Delta\left(t\right) & = & \left[\Delta A\left(t\right)-\Delta B\left(t\right)K\left(t_{s}\right)C\right]^{T}P\nonumber \\
 &  & +P\left[\Delta A\left(t\right)-\Delta B\left(t\right)K\left(t_{s}\right)C\right]
\end{eqnarray}

\noindent where, $\Delta A\left(t\right)=A\left(t\right)-A\left(t_{s}\right)$
and $\Delta B\left(t\right)=B\left(t\right)-B\left(t_{s}\right)$.
Taking norm on both side of equation (23) yields

\vspace{-1.2em}

\begin{eqnarray}
\begin{Vmatrix}\Delta\left(t\right)\end{Vmatrix} & \leq & 2\begin{Vmatrix}P\left[\Delta A\left(t\right)-\Delta B\left(t\right)K\left(t_{s}\right)C\right]\end{Vmatrix}\nonumber \\
 & \leq & 2\begin{Vmatrix}P\end{Vmatrix}\left(\begin{Vmatrix}\Delta A\left(t\right)\end{Vmatrix}+\begin{Vmatrix}\Delta B\left(t\right)\end{Vmatrix}\begin{Vmatrix}K\left(t_{s}\right)\end{Vmatrix}\begin{Vmatrix}C\end{Vmatrix}\right)\nonumber \\
 & \leq & 2\lambda_{M}\left(P\right)\delta\Delta t
\end{eqnarray}

\noindent where, $\delta=\delta_{A}+\delta_{B}k_{M}$ and $\Delta t=t-t_{s}$.
Note that $\begin{Vmatrix}C\end{Vmatrix}=1$. Also $\left\Vert P\right\Vert =\lambda_{M}\left(P\right)$
for all $P\in\mathbf{S}^{+}$. Getting back to equation (22) we have

\vspace{-1.2em}

\begin{eqnarray}
\dot{E} & = & \dot{E}\left(t_{s}\right)+x^{T}\Delta\left(t\right)x\leq-s\begin{Vmatrix}x\end{Vmatrix}^{2}+\begin{Vmatrix}\Delta\left(t\right)\end{Vmatrix}\begin{Vmatrix}x\end{Vmatrix}^{2}\nonumber \\
 & \leq & -s\begin{Vmatrix}x\end{Vmatrix}^{2}+\begin{Vmatrix}\Delta\left(t\right)\end{Vmatrix}\begin{Vmatrix}x\end{Vmatrix}^{2}\nonumber \\
 & \leq & -\left[s-2\lambda_{M}\left(P\right)\delta\Delta t\right]\begin{Vmatrix}x\end{Vmatrix}^{2}\nonumber \\
 & \leq & -\left(\alpha-2\delta\Delta t\right)E
\end{eqnarray}

To find the minimum possible value of $\tau$ we can substitute $\dot{E}=-\gamma\alpha E$
in inequality (25). This gives,

\noindent 
\begin{eqnarray}
\tau & \geq & \frac{\alpha\left(1-\gamma\right)}{2\delta}
\end{eqnarray}

Hence, $\tau_{min}=\frac{\alpha\left(1-\gamma\right)}{2\delta}$.
Let $E=E_{r}$ at the end of rest mode. We integrate inequality (25)
from $t=t_{s}$ to $t=t_{s}+\frac{\alpha\left(1-\gamma\right)}{2\delta}$
to find the minimum possible contraction in rest mode,

\noindent \vspace{-2.3em}

\begin{eqnarray}
E_{r} & \leq & \beta_{r}E_{s}\quad\text{{where,}}
\end{eqnarray}

\noindent $\beta_{r}=\exp\left(-\frac{\alpha^{2}\left(1-\gamma^{2}\right)}{4\delta}\right)$,
the worst case contraction factor. Using inequality (20) and (27)
we get
\begin{eqnarray}
E_{r} & \leq & \beta E^{*}\quad\text{{where,}}
\end{eqnarray}

\noindent $\beta=\beta_{s}\beta_{r}=\exp\left(-\frac{\alpha^{2}\left(1-\gamma^{2}\right)}{4\delta}+\frac{2\lambda_{M}^{\mathcal{L}}T}{\lambda_{m}\left(P\right)}\right)$.
If $\beta\in\left(0,1\right)$, there will be an overall decrease
in $E$ in a rest-scan cycle. $\beta\in\left(0,1\right)$ can be assured
if
\begin{equation}
T<\frac{\lambda_{m}\left(P\right)\alpha^{2}\left(1-\gamma^{2}\right)}{8\delta\lambda_{M}^{\mathcal{L}}}
\end{equation}

\noindent Now we want to prove that inequality (29) ensures that $E\left(t\right)\rightarrow0$
(and hence%
\footnote{\noindent As $E=x^{T}Px\geq\lambda_{m}\left(P\right)\left\Vert x\right\Vert ^{2}$
it implies $\left\Vert x\right\Vert \leq\sqrt{\frac{E}{\lambda_{m}\left(P\right)}}$
where $\lambda_{m}\left(P\right)>0$ as $P\succ0$. Hence if $E\rightarrow0$
then $\left\Vert x\right\Vert \rightarrow0$ implying $x\rightarrow0$.%
} $x\left(t\right)\rightarrow0$) as $t\rightarrow\infty$.Before proceeding
forward we note two things: 1) Theoretically predictable duration
of a rest-scan cycle is $T_{rs}=2T+\frac{\alpha\left(1-\gamma\right)}{2\delta}$.
2) $T_{rs}$ is a conservative estimate of the duration of a rest-scan
cycle. Hence it may as well happen that the rest mode lasts for a
longer time leading to unpredictable contraction. This is clearly
shown in Fig. 7. Now we want to put an upper bound on $E\left(t\right)$.
Say we want to upper bound $E\left(t\right)$ in any one blue dots
shown in Fig. 7, i.e. in the green zones. This can be done as follows

\noindent \vspace{-1.5em}

\begin{equation}
E\left(t\right)\leq E_{o}\beta^{\eta\left(t\right)}\exp\left(-\gamma\alpha\left(t-\eta\left(t\right)T_{rs}\right)\right)
\end{equation}

\noindent where, $E=E_{o}$ at $t=0$, $\eta\left(t\right)$ is the
number of \textit{complete} rest-scan cycle%
\footnote{Without additional knowledge of system dynamics, $\eta\left(t\right)$
is not predictable.%
} before time $t$. In inequality (30), $\beta^{\eta\left(t\right)}$
is the predictable contraction factor contributed by the red zones
in Fig. 7 and the last term is the unpredictable contraction factor
contributed by the green zones in Fig. 7. Note that in green zones
all we can assure is that, $\dot{E}\leq-\gamma\alpha E$, which when
integrated yields the last term of inequality (30). Now we want to
upper bound $E\left(t\right)$ in one of the red dots shown in Fig.
7 which can be done as

\noindent \vspace{-2.0em}

\begin{equation}
E\left(t\right)\leq E_{o}\beta_{s}\beta^{\eta\left(t\right)}\exp\left(-\gamma\alpha\left(t-T_{rs}-\eta\left(t\right)T_{rs}\right)^{+}\right)
\end{equation}

In inequality (31), the operator $\left(a\right)^{+}=\max\left(0,\, a\right)$.
Inequality (31) resembles (30) except that in this case the system
may be in scan mode leading to the extra expansion factor $\beta_{s}$
in the expression. There is an extra $-T_{rs}$ in the last term to
deduct the predictable time of the current rest-scan cycle. Among
inequality (30) and (31) , (31) is definitely the worst case upper
bound on $E\left(t\right)$ due to the presence of two additional
terms, i.e. $-T_{rs}$ inside $\exp\left(-\gamma\alpha\left(\cdot\right)\right)$
and $\beta_{s}\geq1$. From inequality (31) it is clear that if $\beta\in\left(0,1\right)$
then, $E\left(t\right)\rightarrow0$ as $t\rightarrow\infty$. This
completes the proof of Theorem 1.

\noindent 
\begin{figure}[t]
\begin{centering}
\includegraphics[width=9cm,height=6cm]{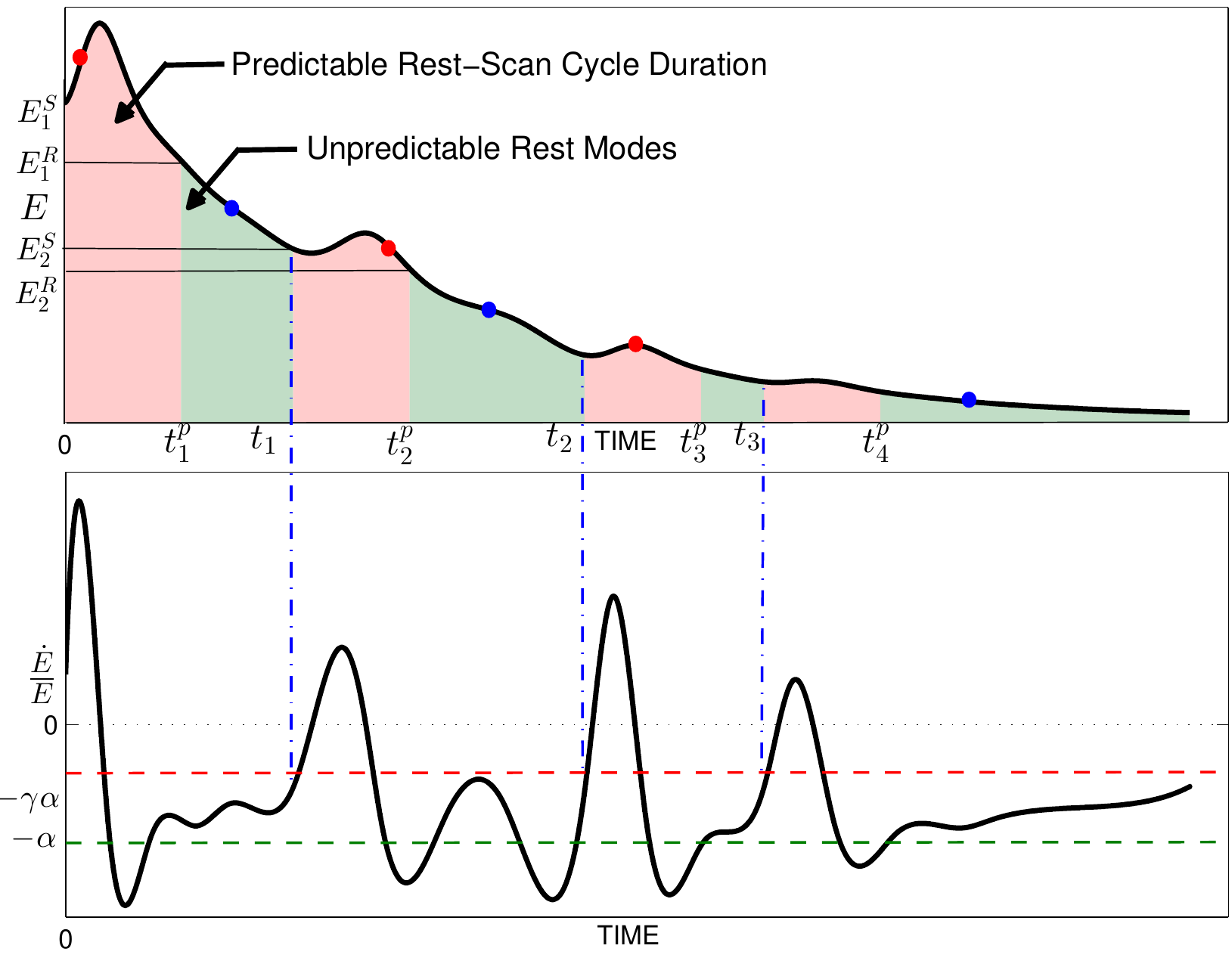}
\par\end{centering}

\caption{Plot showing the actual and predictable duration of rest-scan cycle.
Red zones shows the predictable duration, i.e. $\left(t_{i}^{p}-t_{i-1}\right)=T_{rs}$,
while the actual duration is $t_{i}-t_{i-1}$. Also note that $E_{i}^{R}\leq\beta E_{i}^{S}$,
the theoretically predictable contraction in a rest-scan cycle. The
green zones shows the unpredictable rest modes.}

\vspace{-1.5em}
\end{figure}

\vspace{-1.2em}

\textit{Remark 4: }From inequality (31) convergence is better if
$\beta$ is small. According to inequality (28), $\beta$ decreases
as $\alpha$ (or $s$) increases and $T$ decreases. The effect of
$\gamma$ on $\beta$ is more involved. In inequality (31), if $\gamma$
increases then $\beta$ increases (refer inequality (28)). Thus predictable
contraction decreases (due to large $\beta$) but unpredictable contraction
increases (due to large $\gamma$). Hence $\gamma$ can be neither
too high nor too low.

\textit{Remark 5: }Note that RGS strategy doesn't impose any theoretical
limit on the rate of variation, $\delta_{A}$ and $\delta_{B}$, of
LTV system. This is perhaps one of the novelty of RGS.

\textit{Remark 6 (RGS for LTI Systems): }For uncertain LTI systems,
\textit{Theorem 1} will have \textbf{[C1]} without any change. However,
\textbf{[C2]} will no longer impose an upper bound on $T$ but will
just demand a finite non-zero value. This will ensure that a stabilizing
gain is found within a finite time period of $2T$.

\subsection{Bilinear Matrix Inequality(BMI) Optimization Problem\label{sub:optimization}}

Foundation of RGS is based on the validity of \textbf{[C1]} for a
given uncertain LTV/LTI system. We pose a BMI optimization problem
to check the validity of \textbf{[C1]} and in the process find the
value of $P$ and $s$ needed for implementing RGS.

We start our discussion by formally defining the set $\mathcal{L}$.
Let $A\left(t\right)$ and $B\left(t\right)$ be functions of $p$
independent physical parameters $\Theta=\left[\Theta_{1},\:\Theta_{2},\ldots,\Theta_{p}\right]^{T}$,
i.e. $\left(A\left(t\right),B\left(t\right)\right)=F\left(\Theta\left(t\right)\right)$.
$\Theta\left(t\right)$ is time varying and at every time instant
is also associated with a uncertainty set because of parameter uncertainty.
We assume that every physical parameter $\Theta_{i}$ is bounded,
i.e. $\Theta_{i}\left(t\right)\in\left[\Theta_{i}^{L},\Theta_{i}^{H}\right]$.
Then $\Theta\left(t\right)\in\mathcal{S}$, where $\mathcal{S}=\left[\Theta_{1}^{L},\Theta_{1}^{H}\right]\times\left[\Theta_{2}^{L},\Theta_{2}^{H}\right]\times\ldots\times\left[\Theta_{p}^{L},\Theta_{p}^{H}\right]$
a $p-$dimensional hypercube. We assume no knowledge of the time variation
of $\Theta$, i.e. $\Theta\left(t\right)$, but we assume the knowledge
of $\mathcal{S}$. Then $\mathcal{L}$ is the image of $\mathcal{S}$
under the transformation $F$, i.e. $F:\mathcal{S}\rightarrow\mathcal{L}$
where $\mathcal{S}\subset\mathbb{R}^{p}$ and $\mathcal{L}\subset\mathbb{R}^{N\times N}\times\mathbb{R}^{N\times1}$.

\textit{Remark 7: }The system described by equation (15) can be represented
using the compact notation $\begin{bmatrix}a_{1} & a_{2} & \cdots & a_{N} & | & b\end{bmatrix}$.
Hence one can assume that $\mathcal{L}\subset\mathbb{R}^{N+1}$ rather
than $\mathcal{L}\subset\mathbb{R}^{N\times N}\times\mathbb{R}^{N\times1}$.
This reduces notational complexity by making the elements of $\mathcal{L}$
a vector rather than ordered pair of two matrix.

At this point we would be interested in formulating \textbf{[C1]}
as an optimization problem. With a slight misuse of variable $s$
we can state the following problem.

\begin{problem}
\textbf{[C1]} holds \textit{if and only if} the optimal value of the problem\\
\textit{minimize: }$s$\\
\textit{subject to:}\\
$\left(A-BK_{AB}C\right)^{T}P+P\left(A-BK_{AB}C\right)\preceq sI\quad\forall\left(A,B\right)\in\mathcal{L}$
$P\succ0$\\
with the design variables $s\in\mathbb{R}$, $P\in\mathbb{R}^{N\times N}$ and $K_{AB}\in\left[k_{m},k_{M}\right]$ is negative.
\end{problem}

Note the use of $K_{AB}$ instead of $K$ in Problem 1. It is to signify
that we do not need a common gain $K$ for all $\left(A,B\right)\in\mathcal{L}$.
Perhaps we can have seperate gains $K_{AB}$ for every $\left(A,B\right)\in\mathcal{L}$
satisfying Problem 1 and the RGS strategy described in Section~\ref{sub:Reflective-Gain-Space}
will search for it. However the optimization problem described in
Problem 1 is semi-infinite%
\footnote{Semi-Infinite optimization problems are the ones with infinite constraints.%
} and is hence not computationally tractable. We will now pose \textit{Problem 1}
as a finite optimization problem.

We can always bound the compact set $\mathcal{L}$ by a convex polytope.
Define a polytope%
\footnote{For a given $\mathcal{L}$, $\mathcal{P}$ is not unique.%
} $\mathcal{P}=Conv\left\{ \mathcal{V}\right\} $ where $\mathcal{V}=\left\{ \left(A_{i},B_{i}\right):i=1,2,\ldots,m\right\} $,
the $m$ vertices of the convex polytope s.t., if $\left(A,B\right)\in\mathcal{L}$
then $\left(A,B\right)\in\mathcal{P}$. Then $\mathcal{L}\subseteq\mathcal{P}$.

We now give the following example to illustrate concepts related to
$\mathcal{S}$, $\mathcal{L}$ and $\mathcal{P}$ (discussed above)
and also discuss how to calculate $\delta_{A}$ and $\delta_{B}$
(discussed in \textit{Assumption 3}).

\noindent 
\begin{figure}[t]
\begin{centering}
\includegraphics[width=9cm,height=5cm]{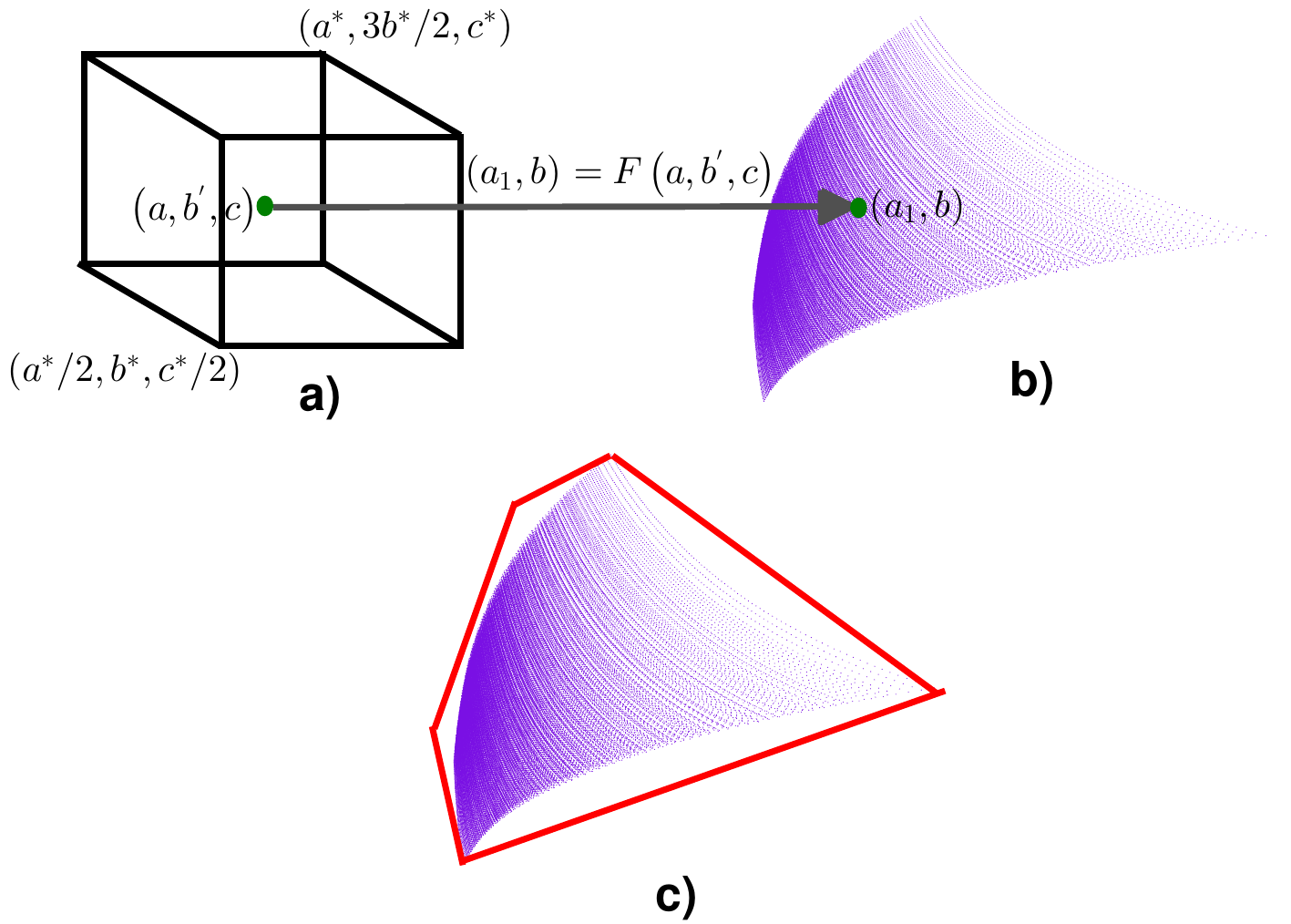}
\par\end{centering}

\caption{a) Uncertainty set $\mathcal{S}$ associated with the physical parameters.
The coordinates of two diagonally opposite vertices is shown in the
figure. b) Uncertainty set $\mathcal{L}$ associated with the LTV
system characterized by equation (15). c) Bounding Convex Polytope
$\mathcal{P}$ of set $\mathcal{L}$. The red lines are the edges
of the polytope $\mathcal{P}$. $\mathcal{P}$ has $5$ vertices which
forms the elements of set $\mathcal{V}$.}

\vspace{-1.2em}
\end{figure}

\vspace{-1.em}

\textit{Example 2: }Consider the following scalar LTV system:

\begin{equation}
\dot{x}_{1}=\frac{a\left(t\right)^{3}c\left(t\right)^{2}}{b^{'}}x_{1}+\frac{\sqrt{b^{'}c\left(t\right)}}{a\left(t\right)^{\nicefrac{2}{3}}}u\:;\quad y=x_{1}
\end{equation}

Here $a$, $b^{'}$, $c$ are the physical parameters which may represent
quantities like mass, friction coefficient, resistance, capacitance
etc. $a$ and $c$ are time varying while $b^{'}$ is an uncertain
parameter. $a$ and $c$ varies as: $a\left(t\right)=a^{*}-\frac{a^{*}}{2}\exp\left(-\frac{t}{\tau_{a}}\right)$,
$c\left(t\right)=\frac{c^{*}}{2}+\frac{c^{*}}{2}\exp\left(-\frac{t}{\tau_{c}}\right)$
and $b^{'}$ lies in the uncertainty set $\left[b^{*},\,\frac{3b^{*}}{2}\right]$.
Here physical parameter $\Theta=\left[a,\, b,\, c\right]$ and hence
$p=3$. From the time variation of $a$ and $c$ we can conclude that
$a\in\left[\frac{a^{*}}{2},\, a^{*}\right]$ and $c\in\left[\frac{c^{*}}{2},\, c^{*}\right]$.
Therefore, $\mathcal{S}\in\left[\frac{a^{*}}{2},\, a^{*}\right]\times\left[b^{*},\,\frac{3b^{*}}{2}\right]\times\left[\frac{c^{*}}{2},\, c^{*}\right]$
which is an hypercube as shown in Fig. 8a.

To find set $\mathcal{L}$ first note that for this example $N=1$
as the system described by equation (32) is scalar. Therefore $\mathcal{L}\subset\mathbb{R}^{2}$
and every element $\left[a_{1},\, b\right]^{T}\in\mathcal{L}$ is
a map $\left[a_{1},\, b\right]^{T}=F\left(a,\, b^{'},\, c\right)$
where

\vspace{-0.4em}

\begin{equation}
F\left(a,\, b^{'},\, c\right)=\begin{bmatrix}\frac{a{}^{3}c^{2}}{b^{'}}\\
\frac{\sqrt{b^{'}c}}{a{}^{\nicefrac{2}{3}}}
\end{bmatrix}
\end{equation}

\vspace{-0.4em}

One method to obtain set $\mathcal{L}$ is to divide the hypercube
$\mathcal{S}$ into uniform grids and then map each grid using equation
(33). Set $\mathcal{L}$ for this example is shown in Fig. 8b. Though
this method to obtain $\mathcal{L}$ from $\mathcal{S}$ is computationally
expensive, one must appreciate that for a general system there is
no other elegant method. Now we want to find a convex polytope $\mathcal{P}$
which bounds $\mathcal{L}$. One such polytope is shown in Fig. 8c.
This polytope has $m=5$ vertices. In practise polytope $\mathcal{P}$
can be found using \textit{convex-hull functions} available widely
in many programming languages. One popular reference is \cite{qhull}.

We now discuss how to calculate $\delta A$ and $\delta B$ for this
example. At this point one must appreciate that in most practical
scenario the controller designer may not explicitly know the equations
governing the rate of change of physical parameters, i.e. say in this
example the designer may not know $a\left(t\right)$ and $c\left(t\right)$
explicitly. However it seems practical to assume knowledge of the
bounds on the rate of change of physical parameters, i.e. controller
designer may know that $0<\dot{a}\leq\frac{a^{*}}{2\tau_{a}}$ and
$-\frac{c^{*}}{2\tau_{c}}\leq\dot{c}<0$. There is no standard way
to calculate $\delta A$ and $\delta B$ from these bounds. Also dependent
on the method used, one may get different estimate of $\delta A$
and $\delta B$. For this example $\delta A$ and $\delta B$ is calculated
as follows.

\vspace{0.5em}

$\delta_{A}=\max\left(\left|\frac{d}{dt}\left(\frac{a\left(t\right)^{3}c\left(t\right)^{2}}{b^{'}}\right)\right|\right)=\max\left(\frac{\left|3c^{2}a^{2}\dot{a}+2a^{3}c\dot{c}\right|}{b^{'}}\right)$

\vspace{0.5em}

$\quad\;=\frac{\max\left(3\sup\left(c\right)^{2}\sup\left(a\right)^{2}\max\left(\dot{a}\right),\:2\sup\left(a\right)^{3}\sup\left(c\right)\max\left(\left|\dot{c}\right|\right)\right)}{\inf\left(b^{'}\right)}$

\vspace{0.5em}

$\quad\;=\frac{\max\left(3\left(c^{*}\right)^{2}\left(a^{*}\right)^{2}\left(\frac{a^{*}}{2\tau_{a}}\right),\:2\left(a^{*}\right)^{3}\left(c^{*}\right)\left(\frac{c^{*}}{2\tau_{c}}\right)\right)}{\left(b^{*}\right)^{2}}$

\vspace{0.5em}

$\quad\;=\frac{\left(a^{*}\right)^{3}\left(c^{*}\right)^{2}}{2\left(b^{*}\right)^{2}}\max\left(\frac{3}{\tau_{a}},\:\frac{2}{\tau_{c}}\right)$

\vspace{0.8em}

$\delta_{B}=\max\left(\left|\frac{d}{dt}\left(\frac{\sqrt{b^{'}c\left(t\right)}}{a\left(t\right)^{\nicefrac{2}{3}}}\right)\right|\right)$

\vspace{0.5em}

$\quad\;=\max\left(\sqrt{b^{'}}\left|-\frac{2\sqrt{c}\dot{a}}{3a^{\nicefrac{5}{3}}}+\frac{\dot{c}}{2\sqrt{c}a^{\nicefrac{2}{3}}}\right|\right)$

\vspace{0.5em}

$\quad\;=\sqrt{\sup\left(b^{'}\right)}\left(\frac{2\sqrt{\sup\left(c\right)}\left(\frac{a^{*}}{2\tau_{a}}\right)}{3\inf\left(a\right)^{\nicefrac{5}{3}}}+\frac{\left(\frac{c^{*}}{2\tau_{c}}\right)}{2\sqrt{\inf\left(c\right)}\inf\left(a\right)^{\nicefrac{2}{3}}}\right)$

\vspace{0.5em}

${\displaystyle \quad\;=\frac{\sqrt{c^{*}}}{\left(a^{*}\right)^{\nicefrac{2}{3}}}}\left(\frac{2^{\nicefrac{5}{3}}}{3\tau_{c}}+\frac{1}{2^{\nicefrac{5}{6}}\tau_{c}}\right)$

\begin{lemma}
Under \textit{Assumption 4}, if for a given $P\in\mathbb{R}^{N\times N}$ and $s\in\mathbb{R}$ there exist a $K_{i}\in\left[k_{m},k_{M}\right]$ s.t.
\[ \left(A_{i}-B_{i}K_{i}C\right)^{T}P+P\left(A_{i}-B_{i}K_{i}C\right)\preceq sI,\; \left(A_{i},B_{i}\right)\in\mathcal{V} \]
then there exist a $K_{AB}\in\left[k_{m},k_{M}\right]$ s.t.
\[ \left(A-BK_{AB}C\right)^{T}P+P\left(A-BK_{AB}C\right)\preceq sI\quad\forall\left(A,B\right)\in\mathcal{L} \]
\end{lemma}

\textit{Proof: }We first define a $N\times N$ matrix $\Gamma$ all
elements of which are $0$ except its $\left(N,1\right)$ element
which is $1$. Also note that all $\left(A,B\right)\in\mathcal{L}$
can be written as a convex combination of the elements $\mathcal{V}$
of the convex polytope $\mathcal{P}$. Mathematically, $\forall\left(A,B\right)\in\mathcal{L}$
there exists scalars $\theta_{i}\geq0,\, i=1,\,2,\ldots m$ s.t.

\vspace{-0.4em}

\[
\left(A,B\right)=\sum_{i=1}^{m}\theta_{i}\left(A_{i},B_{i}\right)\;\text{{and}}\;\sum_{i=1}^{m}\theta_{i}=1
\]

\noindent where, $\left(A_{i},\, B_{i}\right)\in\mathcal{V}\;\forall i=1,\,2,\ldots,\, m$.
Now,\\
\vspace{-0.2em}

$\left(A-BK_{AB}C\right)^{T}P+P\left(A-BK_{AB}C\right)$\\
\vspace{-0.2em}

$=\left(A^{T}P+PA\right)-bK_{AB}\left(\Gamma^{T}P+P\Gamma\right)$\\
\vspace{-0.5em}
\begin{equation}
=\sum_{i=1}^{m}\theta_{i}\left(A_{i}^{T}P+PA_{i}\right)-{\displaystyle \left(\Gamma^{T}P+P\Gamma\right)K_{AB}}\sum_{i=1}^{m}\theta_{i}b_{i}
\end{equation}

\noindent Equation (34) is possible because $BC=b\left(t\right)\Gamma$
owing to the special structure of $B$ and $C$ matrix. Using the
inequality\\
\vspace{-0.4em}

$\left(A_{i}-B_{i}K_{i}C\right)^{T}P+P\left(A_{i}-B_{i}K_{i}C\right)\preceq sI,\;\left(A_{i},B_{i}\right)\in\mathcal{V}$\\
\vspace{-0.4em}

\noindent in equation (34) we have\\
\vspace{-0.6em}

$\left(A-BK_{AB}C\right)^{T}P+P\left(A-BK_{AB}C\right)$

\vspace{-0.6em}

\begin{equation}
\preceq\left[\sum_{i=1}^{m}\theta_{i}K_{i}b_{i}-K_{AB}\sum_{i=1}^{m}\theta_{i}b_{i}\right]\left(\Gamma^{T}P+P\Gamma\right)+sI
\end{equation}

\noindent Choosing, $K_{AB}=\frac{{\displaystyle \sum_{i=1}^{m}}\theta_{i}K_{i}b_{i}}{{\displaystyle \sum_{i=1}^{m}}\theta_{i}b_{i}}$
in inequality (35) yields\\
\\
\vspace{-0.4em}

$\quad\quad\quad\left(A-BK_{AB}C\right)^{T}P+P\left(A-BK_{AB}C\right)\preceq sI$\\
\vspace{-0.2em}

Now all we need to do is to prove that the chosen $K_{AB}$ lies in
the interval $\left[k_{m},k_{M}\right]$. We know that $k_{m}\leq K_{i}\leq k_{M}$.
Therefore under \textit{Assumption 4}%
\footnote{Without Assumption 4, the denominator $\sum_{i=1}^{m}\theta_{i}b_{i}$
may become $0$.%
},\\
\vspace{-0.4em}

$\frac{{\displaystyle \sum_{i=1}^{m}\theta_{i}k_{m}b_{i}}}{{\displaystyle \sum_{i=1}^{m}\theta_{i}b_{i}}}\leq K_{AB}\leq\frac{{\displaystyle \sum_{i=1}^{m}\theta_{i}k_{M}b_{i}}}{{\displaystyle \sum_{i=1}^{m}\theta_{i}b_{i}}}\Rightarrow k_{m}\leq K_{AB}\leq k_{M}$\\
\\
\vspace{-0.6em}

This concludes the proof of Lemma 3. The importance of Lemma 3 is
that it reduces the semi-infinite optimization problem posed in Problem
1 into a finite optimization problem. This results into the most important
result of this section.

\begin{theorem}
\textbf{[C1]} holds \textit{if} the optimal value of the problem\\
\textit{minimize: }$s$\\
\textit{subject to:}\\
$\left(A_{i}-B_{i}K_{i}C\right)^{T}P+P\left(A_{i}-B_{i}K_{i}C\right)\preceq sI\quad\forall\left(A_{i},B_{i}\right)\in\mathcal{V}$
$k_{m}\leq K_{i}\leq k_{M},\quad i=1,2,\ldots,m$\\
$P\succ0$\\
with design variables $s\in\mathbb{R}$, $P\in\mathbb{R}^{N\times N}$ and $K_{i}\in\mathbb{R},\: i=1,2,\ldots,m$ is negative.
\end{theorem}

\textit{Proof: }The proof follows from \textit{Problem 1} and \textit{Lemma 3}.

Note that while Problem 1 is a \textit{necessary and sufficient}
condition, Theorem 2 is a \textit{sufficient} condition. This is
due to the fact that $\mathcal{L}\subseteq\mathcal{P}$ leading to
some conservativeness in the optimization problem proposed in Theorem
2.

Theorem 2 poses the classical \textit{BMI Eigenvalue Minimization Problem (BMI-EMP)}
in variables $P$ and $K_{i}$. As such BMI's are non-convex in nature
leading to multiple local minimas. Several algorithms to find the
local minima exist in literature (see \cite{localBMI99_pathfollowing,localBMI98_ILMI}).
Algorithms for global minimization of BMI's is rather rare and have
received attention in works like \cite{goh1994global,globalBMI97_fujioka}.
Our approach is similar to \cite{goh1994global} which is basically
a \textit{Branch and Bound} algorithm. Such an algorithm works by
bounding $s$ by a lower bound $\Phi_{L}$ and an upper bound $\Phi_{U}$,
i.e. $\Phi_{L}\leq s\leq\Phi_{U}$. The algorithm then progressively
refines the search to reduce $\Phi_{U}-\Phi_{L}$. Our main algorithm
consists of \textit{Algorithm 4.1} of \cite{goh1994global} (Page
4). The \textit{Alternating SDP method} mentioned in \cite{globalBMI07_parallel}
(Page 2) is used for calculating $\Phi_{U}$. For calculating $\Phi_{L}$
we have used the \textit{convex relaxation} technique first introduced
in \cite{globalBMI97_fujioka} (also discussed in \cite{globalBMI07_parallel},
equation (9)). In Appendix \ref{sec:BMIEMP} we present a working
knowledge of our algorithm. For detailed explanation the readers may
refer the corresponding work \cite{goh1994global,globalBMI07_parallel}.

\textit{Theorem 2} poses an optimization problem with $\left(A_{i},B_{i}\right)\in\mathcal{V}$,
$k_{m}$ and $k_{M}$ as inputs and $P$, $s$ and $K_{i},\, i=1,2,\ldots,m$
as outputs. But $k_{m}$ and $k_{M}$ are not known. An initial estimate
of $k_{m}=k_{m}^{RH}$ and $k_{M}=k_{M}^{RH}$ is obtained by using
Routh-Hurwitz criteria%
\footnote{Routh-Hurwitz criteria is used to find the bounds on the feedback
gains for which a SISO LTI system is closed loop stable. Refer \cite{ogata}
for details.%
} for each $\left(A_{i},B_{i}\right)\in\mathcal{V}$. Let $K_{i}=K_{i}^{RH},\, i=1,2,\ldots,m$
be the output of the optimization problem with this initial estimate.
Let $k_{m}^{*}=\underset{1\leq i\leq m}{\min}\left(K_{i}^{RH}\right)$
and $k_{M}^{*}=\underset{1\leq i\leq m}{\max}\left(K_{i}^{RH}\right)$,
then the following holds:\\
\vspace{-0.2em}

\noindent 1) $\ensuremath{k_{m}^{*}\geq k_{m}^{RH}}$ and $\ensuremath{k_{M}^{*}\leq k_{M}^{RH}}$.
This is because $\ensuremath{k_{m}^{RH}}\leq K_{i}^{RH}\leq\ensuremath{k_{M}^{RH}}\,,\;\forall\, i=1,\,2,\,\ldots,\, m$
and hence\\
\vspace{-0.3em}

\noindent $\qquad\quad k_{m}^{RH}\leq\underset{1\leq i\leq m}{\min}\left(K_{i}^{RH}\right)\leq\underset{1\leq i\leq m}{\max}\left(K_{i}^{RH}\right)\leq k_{M}^{RH}$\\
\vspace{-0.3em}

\noindent 2) The outputs, $P$, $s$ and $K_{i}$, obtained by running
the optimization algorithm with \textit{a)} $k_{m}=k_{m}^{RH}$ and
$k_{M}=k_{M}^{RH}$ or \textit{b)} $k_{m}=k_{m}^{*}$ and $k_{M}=k_{M}^{*}$,
will be the same. This is because the gains $K_{i}^{RH}$ obtained
by running the optimization algorithm with $k_{m}=k_{m}^{RH}$ and
$k_{M}=k_{M}^{RH}$ also satisfies the bounds $k_{m}^{*}\leq K_{i}^{RH}\leq k_{M}^{*}$.\\
\vspace{-0.6em}

Therefore if we choose $k_{m}=k_{m}^{*}$ and $k_{M}=k_{M}^{*}$ we
would get a smaller RGS gain set, i.e. $\left(k_{M}^{*}-k_{m}^{*}\right)\leq\left(k_{M}^{RH}-k_{m}^{RH}\right)$,
without compromising the convergence coefficient $s$. A smaller RGS
gain set will ease the controller design in an analog setting.

We now give a bound on $\lambda_{M}^{\mathcal{L}}$ (defined in \textit{Theorem 1}).
It is not possible to calculate $\lambda_{M}^{\mathcal{L}}$ with
the formula given in \textit{Theorem 1} as it will involve search
over the dense set $S_{Q}^{\mathcal{L}}$. Define a set

$\quad S_{Q}^{\mathcal{P}}=\left\{ Q\,:\: Q=\left(A-BKC\right)^{T}P+P\left(A-BKC\right)\right.$

$\quad\quad\quad\;\quad\left.\vphantom{A_{C}\left(t\right)^{T}}\left(A,B\right)\in\mathcal{P},\: K\in\left[k_{m},k_{M}\right]\right\} $

\noindent Let $\lambda_{M}^{\mathcal{P}}=\underset{Q\in S_{Q}^{\mathcal{P}}}{\max}\lambda_{M}\left(Q\right)$.
As $\mathcal{L}\subseteq\mathcal{P}$ it is obvious that, $S_{Q}^{\mathcal{L}}\subseteq S_{Q}^{\mathcal{P}}$
($S_{Q}^{\mathcal{L}}$ defined in Theorem 1). Therefore%
\footnote{This is more like saying that the maximum of a function ($\lambda_{M}\left(\cdot\right)$
here) over a bigger set ($S_{Q}^{\mathcal{P}}$ here) will be greater
than the maximum over a smaller set ($S_{Q}^{\mathcal{L}}$ here).%
}, $\lambda_{M}^{\mathcal{L}}\leq\lambda_{M}^{\mathcal{P}}$. Thus
$\lambda_{M}^{\mathcal{P}}$ gives an estimate of $\lambda_{M}^{\mathcal{L}}$
by upper bounding it. It can be shown that for a scalar gain $K$
and the specific structure of $B$ and $C$ (refer equation (15)),
it can be proved that $S_{Q}^{\mathcal{P}}$ is compact convex set.
Also $\lambda_{M}\left(Q\right)$ is a convex function for all $Q\in\mathbf{S}^{+}$
(refer Page 82 of \cite{boyd2009convex}). It is well known that global
maxima of a convex function over a convex compact set only occurs
at some extreme points of the set (refer \cite{rockafellar1997convex}).
Thus the problem of maximizing $\lambda_{M}\left(\cdot\right)$ over
$Q\in S_{Q}^{\mathcal{P}}$ reduces to maximizing $\lambda_{M}\left(\cdot\right)$
over $Q\in S_{Q}^{\mathcal{V}}$ where

$\quad S_{Q}^{\mathcal{V}}=\left\{ Q\,:\: Q=\left(A-BKC\right)^{T}P+P\left(A-BKC\right)\right.$

$\quad\quad\quad\;\quad\left.\vphantom{A_{C}\left(t\right)^{T}}\left(A,B\right)\in\mathcal{V},\: K\in\left\{ k_{m},k_{M}\right\} \right\} $

\noindent the set of vertices of $S_{Q}^{\mathcal{P}}$. This leads
to the following formula

\noindent \vspace{-1.2em}

\begin{equation}
\lambda_{M}^{\mathcal{L}}\leq\lambda_{M}^{\mathcal{P}}=\underset{Q\in S_{Q}^{\mathcal{\mathcal{V}}}}{\max}\lambda_{M}\left(Q\right)
\end{equation}

Inequality (36) can be used to obtain an estimate of $\lambda_{M}^{\mathcal{L}}$.

\textit{Remark 8: }As mentioned in the beginning of Section \ref{sec:cntrl},
for a controller to be implementable in analog framework it has to
be simple. Though the synthesis of RGS parameters ($k_{m}$, $k_{M}$,
$T$, $\alpha$, $\gamma$ and $P$) is complex, one must appreciate
that designing a controller is a \textit{'one time deal'}. RGS in
itself is a simple gain-scheduled controller governed by equation
(16).

\section{Example\label{sec:Example}}

Parallel Plate Electrostatic Actuator (PPA) shown in Fig. 9a forms
a vital component of several miniaturized systems. We perform regulatory
control of PPA to show effectiveness of the proposed analog architecture
and RGS strategy. PPA's as described in \cite{shapiro2011feedback}
(page 183) follows the following dynamics

\noindent \vspace{-1.4em}

\begin{equation}
m\ddot{y}+b\dot{y}+ky=\frac{\varepsilon A}{2\left(G-y\right)^{2}}V_{s}^{2}
\end{equation}

\noindent which is nonlinear in nature. Plant parameter includes spring
constant $\kappa$, damping coefficient $b$, moving plate mass and
area $m$ and $A$ respectively, permittivity $\varepsilon$ and maximum
plate gap $G$. As we are interested in only regulatory control, we
use the following linearized model

\noindent \vspace{-2.2em}

\begin{equation}
\begin{bmatrix}\dot{x_{1}}\\
\dot{x_{2}}
\end{bmatrix}=\begin{bmatrix}0 & 1\\
-\frac{\kappa\left(G-3G_{o}\right)}{m\left(G-G_{o}\right)} & -\frac{b}{m}
\end{bmatrix}\begin{bmatrix}x_{1}\\
x_{2}
\end{bmatrix}+\begin{bmatrix}0\\
\frac{\sqrt{2\varepsilon A\kappa G_{o}}}{m\left(G-G_{o}\right)}
\end{bmatrix}u_{r}
\end{equation}

\noindent where, $x_{1}$ is the displacement from the operating point
$G_{o}$. Plant output is the moving plate position $y=G_{o}+x_{1}$.
Plant input is $V_{s}=V_{b}+V_{u}$. Comparing with RGS theory, $V_{b}=u_{b}$,
the bias voltage to maintain $G_{o}$ as the operating point and $V_{u}=u_{r}=KV_{e}=K\left(G_{o}-y\right)=-Kx_{1}$,
the regulation voltage supplied by the Memristive AGC. Plant parameter
includes spring constant $\kappa$, damping coefficient $b$, moving
plate mass and area $m$ and $A$ respectively, permittivity $\varepsilon$
and maximum plate gap $G$. For $G_{o}>\frac{G}{3}$, the system has
an unstable pole. We perform regulation around $G_{o}=\frac{2G}{3}$.

The \textit{true} plant parameters are, $m=3\times10^{-3}\, kg$,
$b=1.79\times10^{-2}\, Nsm^{-1}$, $G=10^{-3}\, m$, $A=1.6\times10^{-3}\, m^{2}$.
$G$ and $A$ are uncertain but lie in the set $G\in\mathcal{S}_{G}=\begin{bmatrix}0.5 & 2.0\end{bmatrix}\times10^{-3}\, m$,
$A\in\mathcal{S}_{A}=\begin{bmatrix}1.2 & 1.8\end{bmatrix}\times10^{-3}\, m^{2}$.
$\varepsilon$ varies due to surrounding condition as $\varepsilon\left(t\right)=5\varepsilon_{o}+1.5\varepsilon_{o}\sin\left(7.854t\right)$
where $\varepsilon_{o}$ is the permittivity of vacuum. Spring loosening
causes $\kappa$ to decrease as $\kappa\left(t\right)=0.08+0.087e^{-0.8t}\, Nm^{-1}$.
Now we will discuss the steps involved in implementing RGS in an analog
framework.

\noindent \textbf{Step 1 }(Identify $\mathcal{S}$): $\mathcal{S}$
is the uncertainty set of physical parameters first defined in Page
8. Define set $\mathcal{S}_{\varepsilon}=\begin{bmatrix}3.5\varepsilon_{o} & 6.5\varepsilon_{o}\end{bmatrix}$
and $\mathcal{S}_{\kappa}=\begin{bmatrix}0.08 & 0.167\end{bmatrix}\, Nm^{-1}$.
Then the ordered pair $\left(G,A,\varepsilon,\kappa\right)\in\mathcal{S}=\mathcal{S}_{G}\times\mathcal{S}_{A}\times\mathcal{S}_{\varepsilon}\times\mathcal{S}_{\kappa}$.
Note that here $p=4$.

\noindent \textbf{Step 2 }(Find $\ensuremath{\mathcal{P}}$): To
do this we numerically map $\mathcal{S}$ to $\mathcal{L}$ (as shown
in \textit{Example 2}) and then use \textit{convhulln} function
of MATLAB to find a convex polytope $\mathcal{P}$ s.t. $\mathcal{L}\subseteq\mathcal{P}$.
In this case $\mathcal{P}$ consist of $m=4$ vertices. We explicitly
don't mention the computed $\mathcal{P}$ for the sake of neatness.

\noindent \textbf{Step 3 }(Compute $P$, $s$, $\alpha$, $k_{m}$,
$k_{M}$): Solving optimization problem proposed in \textit{Theorem 2}
we get, $s=0.917$, $k_{m}=8600$ and $k_{M}=86000$ and $P=\begin{bmatrix}0.9937 & 0.0757\\
0.0757 & 0.0895
\end{bmatrix}$. Here, $\lambda_{M}\left(P\right)=1$, hence $\alpha=\frac{s}{\lambda_{M}\left(P\right)}=0.917$.

\begin{figure}[t]
\begin{centering}
\includegraphics[scale=0.44]{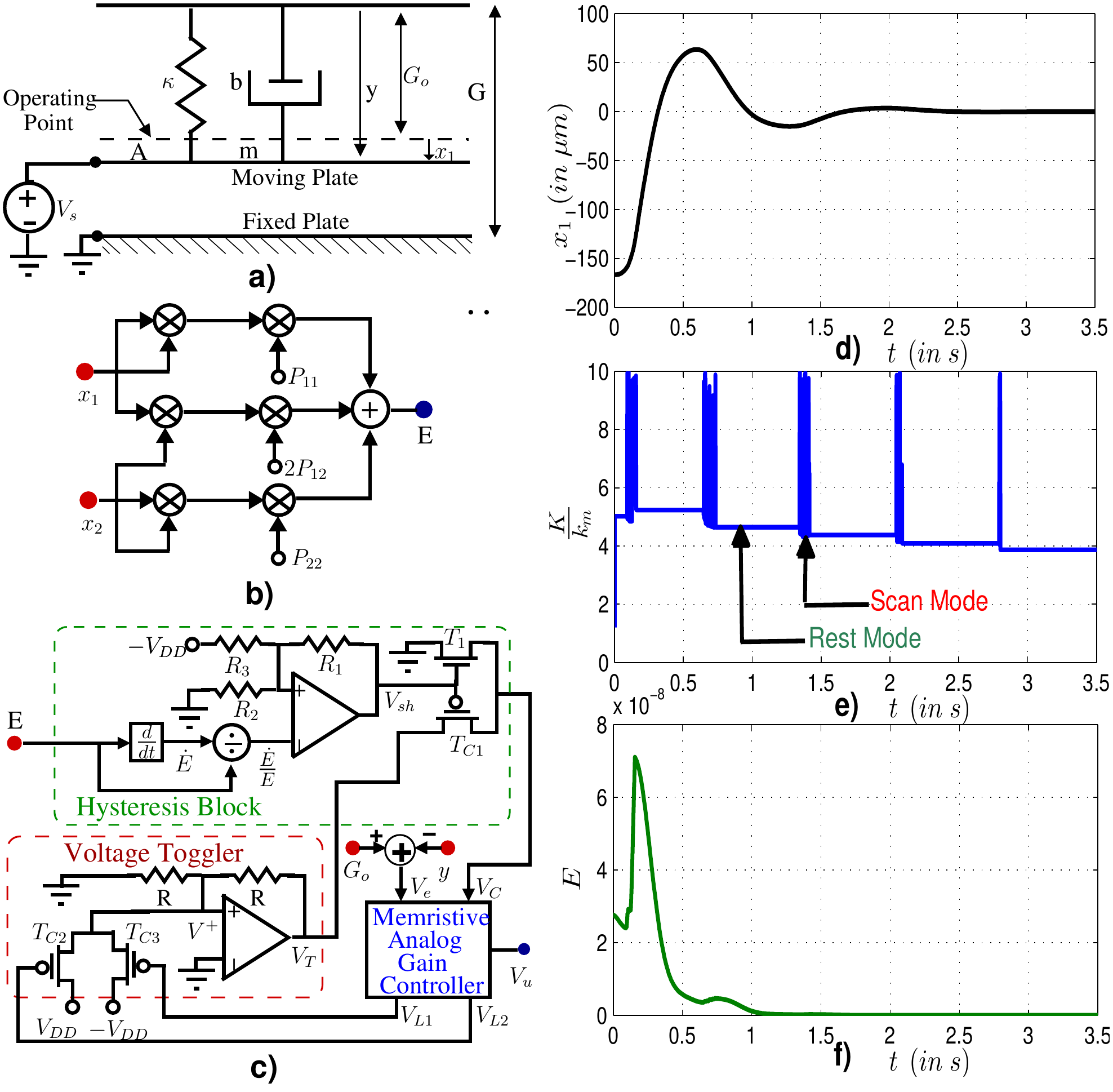}
\par\end{centering}

\caption{a) Schematic of PPA. b) Analog Implementation of $E=x^{T}Px$ for $P\in\mathbb{R}^{2\times2}$. Note that, $P_{12}=P_{21}$ as $P\in\mathbf{S}$. c) Analog Implementation of RGS. d), e), f) Plots of plate position error $x_{1}$ (from operating point $G_{o}$), normalized RGS gain $\frac{K}{k_{m}}$ and Lyapunov Function $E=x^{T}Px$ with respect to time.}

\vspace{-1.5em}
\end{figure}

\noindent \textbf{Step 4 }(Compute $\ensuremath{\lambda_{m}\left(P\right)}$,
$\ensuremath{\lambda_{M}^{\mathcal{L}}}$, $\delta$): For the calculated
$P$, $\lambda_{m}\left(P\right)=0.083$. From equation (36), $\lambda_{M}^{\mathcal{L}}=29.1$.
We will now calculate $\delta_{A}$ and $\delta_{B}$ for this example.
As mentioned in \textit{Example 2} the controller designer does not
know $\kappa\left(t\right)$ and $\varepsilon\left(t\right)$ explicitly
but they do know the bounds on $\dot{\kappa}$ and $\dot{\varepsilon}$
which for this example is: $\left(-0.087\times0.8\right)\leq\dot{\kappa}<0$
and $-\left(1.5\times7.854\times\varepsilon_{o}\right)\leq\dot{\varepsilon}\leq\left(1.5\times7.854\times\varepsilon_{o}\right)$. 

\noindent First note that $\left(2\times2\right)$ element of the
system matrix of the linearized PPA model (described by equation (38))
is not time varying. Hence $\delta_{A}$ can be simply written as

\vspace{0.8em}

$\delta_{A}=\max\left(\left|\frac{d}{dt}\left(-\frac{\kappa\left(t\right)\left(G-3G_{o}\right)}{m\left(G-G_{o}\right)}\right)\right|\right)$

\vspace{0.8em}

$\quad\;=\max\left(\left|\frac{d}{dt}\left(-\frac{\kappa\left(t\right)\left(G-3\left(\frac{2G}{3}\right)\right)}{m\left(G-\frac{2G}{3}\right)}\right)\right|\right)$

\vspace{0.8em}

$\quad\;=\frac{3\max\left(\left|\dot{\kappa}\right|\right)}{m}$

\vspace{0.8em}

$\delta_{B}=\max\left(\left|\frac{d}{dt}\left(\frac{\sqrt{2\varepsilon\left(t\right)A\kappa\left(t\right)G_{o}}}{m\left(G-G_{o}\right)}\right)\right|\right)$

\vspace{0.5em}

$\quad\;=\max\left(\left|\frac{d}{dt}\left(\frac{\sqrt{2\varepsilon\left(t\right)A\kappa\left(t\right)\left(\frac{2G}{3}\right)}}{m\left(G-\frac{2G}{3}\right)}\right)\right|\right)$

\vspace{0.5em}

$\quad\;=\frac{1}{m}\sqrt{\frac{12\sup\left(\mathcal{S}_{A}\right)}{\inf\left(\mathcal{S}_{G}\right)}}\max\left(\left|\dot{\sqrt{\varepsilon K}}\right|\right)$

\vspace{0.5em}

$\quad\;=\frac{1}{m}\sqrt{\frac{12\sup\left(\mathcal{S}_{A}\right)}{\inf\left(\mathcal{S}_{G}\right)}}\max\left(\left|\frac{\kappa\dot{\varepsilon}+\varepsilon\dot{\kappa}}{2\sqrt{\varepsilon}\sqrt{\kappa}}\right|\right)$

\vspace{0.5em}

$\quad\;=\frac{1}{m}\sqrt{\frac{12\sup\left(\mathcal{S}_{A}\right)}{\inf\left(\mathcal{S}_{G}\right)}}\left(\frac{\sup\left(S_{\kappa}\right)\max\left(\left|\dot{\varepsilon}\right|\right)+\sup\left(S_{\varepsilon}\right)\max\left(\left|\dot{\kappa}\right|\right)}{2\sqrt{\inf\left(S_{\varepsilon}\right)}\sqrt{\inf\left(S_{\kappa}\right)}}\right)$

\vspace{0.8em}

\noindent Substituting the values in the above equation of $\delta_{A}$
and $\delta_{B}$ we get $\delta=\delta_{A}+k_{M}\delta_{B}=1351$.

\noindent \textbf{Step 5 }(Compute $T$ and $\gamma$): We arbitarily
choose $\gamma=0.5$. Substituting $\gamma=0.5$ in inequality (29)
we get $T<1.66\times10^{-7}\, s$. Hence we choose $T=10^{-7}\, s$.

\noindent \textbf{Step 6 }(Analog Design): Fig. 9b and 9c combined
shows the analog implementation of RGS. The error voltage $V_{e}$
is the plate position error, i.e. $V_{e}=-x_{1}=\left(G_{o}-y\right)$.
The gain control voltage $V_{C}$ is controlled by the Hysteresis
Block and the Voltage Toggler block. In Rest Mode $V_{C}=0$ thereby
ensuring that the gain $\mathcal{M}$ is constant (refer equation
(12)). In Scan Mode $V_{C}\in\left\{ -V_{DD}\,,\, V_{DD}\right\} $,
$V_{C}=V_{DD}$ to scan from $R_{off}^{S}$ to $R_{on}^{S}$ and $V_{C}=-V_{DD}$
to scan from $R_{on}^{S}$ to $R_{off}^{S}$ (refer equation (12)).
$R_{C}$ controls $\dot{\mathcal{M}}$, the rate of change of gain
(refer equation (12)). Setting $R_{C}=\frac{V_{DD}T}{Q_{M}^{S}}$
ensures a scan time of $T$. The derivation of $R_{C}$ is simple.
Note that the resistance of memristor will change from $R_{off}^{S}$
to $R_{on}^{S}$ if we pass a charge $Q_{M}^{S}$ through it. We want
this change to happen in time $T$ and hence the desired current is
$\frac{Q_{M}^{S}}{T}$. As $V_{C}=V_{DD}$ in scan mode, the required
resistance $R_{C}=\frac{V_{DD}}{\left(\nicefrac{Q_{M}^{S}}{T}\right)}=\frac{V_{DD}T}{Q_{M}^{S}}$.

The Hysteresis Block is a conventional inverting schmitt trigger.
Tuning $R_{1}=\frac{2\left(V_{DD}-\alpha\right)}{\alpha\left(1-\gamma\right)}R_{2}$
and $R_{3}=\frac{2\left(V_{DD}-\alpha\right)}{\alpha\left(1+\gamma\right)}R_{2}$
ensures%
\footnote{The designer is free to choose the resistance $R_{2}$.%
} that the schmitt trigger's output goes from $V_{DD}$ to $-V_{DD}$
at $\frac{\dot{E}}{E}=-\gamma\alpha$ and from $-V_{DD}$ to $V_{DD}$
at $\frac{\dot{E}}{E}=-\alpha$. Due to the transistor arrangement
in Hysteresis Block: 1) $V_{C}=0$ when $V_{sh}=V_{DD}$. Therefore
$V_{sh}=V_{DD}$ implies Rest Mode. 2) $V_{C}=V_{T}$ when $V_{sh}=-V_{DD}$.
It will be explained next that $V_{T}\in\left\{ -V_{DD}\,,\, V_{DD}\right\} $.
Therefore $V_{sh}=-V_{DD}$ implies Scan Mode. So we can conclude
that the Hysteresis Block goes from Rest Mode to Scan Mode when $\frac{\dot{E}}{E}=-\gamma\alpha$
and from Scan Mode to Rest Mode when $\frac{\dot{E}}{E}=-\alpha$.
This is in accordance with the hysteresis shown in Fig. 6a.

$V_{T}$, the output of the Voltage Toggler block, toggles between
$-V_{DD}$ and $V_{DD}$. Recalling equation (12), this will result
in the gain of the memristive gain control block reflect between $\left[\frac{R_{on}^{S}}{R_{I}}\,,\,\frac{R_{off}^{S}}{R_{I}}\right]$.
We now explain the working of this block. Say $V_{T}=V_{DD}$ and
the zone indicating voltages (refer Fig. 3) $V_{L1}=V_{L2}=V_{DD}$.
As $V_{L1}=V_{L2}=V_{DD}$, transistors $T_{C2}$ and $T_{C3}$ are
$OFF$. The voltage $V^{+}=\frac{V_{T}}{2}=\frac{V_{DD}}{2}$. Since
$V^{+}=\frac{V_{DD}}{2}>0$, $V_{T}=V_{DD}$. This shows that the
output of Voltage Toggler block is stable. As $V_{T}=V_{DD}$, memristor's
resistance $M$ will decrease till $M=R_{on}^{S}$ (see equation (12))
at which point $V_{L1}=-V_{DD}$ and $V_{L2}=V_{DD}$ (refer \textit{Case 2}
of \ref{sub:Charge-Saturator-Block}). $T_{C3}$ will be momentarily
$ON$ making $V^{+}=-V_{DD}$ and hence driving $V_{T}$ to $-V_{DD}$.
Then $M$ will increase from $R_{on}^{S}$, making $V_{L1}=V_{L2}=V_{DD}$
and hence driving $T_{C3}$ to $OFF$ state. When this occurs $V^{+}=\frac{V_{T}}{2}=-\frac{V_{DD}}{2}$.
As $V^{+}=-\frac{V_{DD}}{2}<0$, $V_{T}=-V_{DD}$. Similar momentary
transition of $T_{C2}$ to $ON$ state will toggle $V_{T}$ from $-V_{DD}$
to $V_{DD}$ when $M=R_{off}^{S}$.

Several plots corresponding the regulatory performance of PPA under
RGS control strategy is shown in Fig. 9d, e, f. It is interesting
to observe that an LTV system can have multiple rest-scan cycle (see
Fig. 9e). This is because for a time varying system a \textit{stabilizing gain}
at a given instant may not be the \textit{stabilizing gain} at a
later instant due to the change in system dynamics. Unlike LTV system,
a LTI system will have only $1$ rest-scan cycle.

\textit{Remark 9: }In this example $T$ is very low which may seem
to challenge analog design. However for all practical purposes it
is not so. For the sake of simulation we choose a fictitious time
variation of $\kappa\left(t\right)$ and $\varepsilon\left(t\right)$
which is quite fast compared to that found in nature. Therefore $\delta_{A}$
and $\delta_{B}$ is high (refer \textbf{Step 4}) resulting in a
high $\delta$ and hence a low scan time $T$ (refer inequality (29).
In practice, time variation of a system caused by ageing effect and
atmospheric variation is a slow process. Hence, $T$ will be much
higher.

\textit{Remark 10: }To control an array of miniaturized devices (say
PPA) one can reduce the circuitry required by identifying the components
of the circuit which can be common for the entire array. For example,
Synchronization Block can be common for the entire array. Synchronization
of each pair of coupled Memristor Gain Block and Integrator (refer
Fig. 3) can be done in a time multiplexed manner, i.e. each pair of
coupled Memristor Gain Block and Integrator is synchronized using
one Synchronization Block. The oscillator shown in the circuit of
Fig. 3 can also be common for the entire array.

\section{Discussion}

To the best of authors knowledge this paper is one of the first attempts
towards understanding the use of memristor in control applications.
A memristive variable gain architecture has been proposed. We then
propose (and prove) RGS control strategy which can be implemented
using this framework. Simplicity of RGS control strategy is demonstrated
using an example. The extension of this work can take two course.

From \textit{Circuit Standpoint} one may try to design an analog
circuit which mimics the circuit shown in Fig. 2 but with a \textit{lesser number of op-amps}.
Since the synthesis of memristor is still an engineering challenge,
one may speculate regarding the use of variable CMOS resistors (refer
\cite{tunableresistor}) to implement the analog gain controller proposed
in Section \ref{sec:memgain}.

Two milestones have to be addressed before RGS is practically feasible:
1) RGS needs information about the states $x$ which is obtained by
differentiating the output $y$ $N-1\: times$. But differentiation
might amplify the noise. 2) RGS relies on the knowledge of $E$ which
is obtained by performing $x^{T}Px$ using analog circuitry. Such
analog implementation will be easier if $P$ is sparse. Hence from
\textit{Control Theoretic Standpoint}, addressing these two issues
will be the first target of the authors. Later point has been addressed
in \cite{sparse98_Hassibi}. Extending RGS to SISO non-linear and
in general MIMO systems would be the next step. It would also be interesting
to explore other simple control strategies (like \cite{barkana2013simple})
which can be implemented in analog framework.

\appendices{}

\vspace{-1.0em}

\section{Proof of Lemma 1: Drifting Nature of Stabilizing Gain Set $\mathcal{K}_{P,s}\left(t\right)$\label{sec:Proof-Lemma1}}

To prove \textit{Lemma 1} we will take the following steps:

\begin{enumerate}
\item Pick a gain $K\in\mathcal{K}_{P,s}\left(t\right)$.
\item Prove that if $\delta t\rightarrow0$ then there exist a $\Delta K\rightarrow0$ s.t. $\left(K+\Delta K\right)\in\mathcal{K}_{P,s}\left(t+\delta t\right)$.
\item As $\Delta K\rightarrow0$, $\left(K+\Delta K\right)\in\mathcal{K}_{P,s}\left(t\right)$. Hence the gain $K+\Delta K$ belongs to both the sets, $\mathcal{K}_{P,s}\left(t\right)$ and $\mathcal{K}_{P,s}\left(t+\delta t\right)$. This implies that $\mathcal{K}_{P,s}\left(t\right)\bigcap\mathcal{K}_{P,s}\left(t+\delta t\right)\neq\emptyset$.
\end{enumerate}

Now we proceed with the proof. We first pick a gain $\ensuremath{K\in\mathcal{K}_{P,s}\left(t\right)}$.
As $\mathcal{K}_{P,s}\left(t\right)\neq\emptyset,\:\forall t\geq0$
such a gain will exist. For a time $t=t^{*}$ there exists a gain
$K^{*}\in\mathcal{K}_{P,s}\left(t^{*}\right)$ and a scalar $s^{*}>s$
s.t. the following equality holds

\vspace{-1.0em}

\begin{equation}
\left(A^{*}-B^{*}K^{*}C\right)^{T}P+P\left(A^{*}-B^{*}K^{*}C\right)=-s^{*}I
\end{equation}

In equation (39) $A^{*}=A\left(t^{*}\right)$ and $B^{*}=B\left(t^{*}\right)$.
Equation (38) directly follows from the very definition of stabilizing
gain set (defined in page 7). Lets say that at time $t=t^{*}+\delta t$
there exist a $K^{*}+\Delta K$ s.t.

\vspace{-0.8em}

\begin{equation}
\begin{array}{l}
\left[\left(A^{*}+\delta A\right)-\left(B^{*}+\delta B\right)\left(K^{*}+\Delta K\right)C\right]^{T}P\\
+P\left[\left(A^{*}+\delta A\right)-\left(B^{*}+\delta B\right)\left(K^{*}+\Delta K\right)C\right]=-s^{*}I
\end{array}
\end{equation}

In equation (40) $\delta A=\dot{A}\left(t^{*}\right)\,\delta t$ and
$\delta B=\dot{B}\left(t^{*}\right)\,\delta t$ are infinitesimal
change in $A$ and $B$ respectively. Note that $\Delta K$ is a scalar.
Hence, substituting equation (39) in (40) yields

\vspace{-1.0em}

\begin{equation}
\Delta K\mathcal{R}=\left(\delta A-\delta BK^{*}C\right)^{T}P+P\left(\delta A-\delta BK^{*}C\right)
\end{equation}

\noindent where, $\mathcal{R}=C^{T}\left(B^{*}+\delta B\right)^{T}P+P\left(B^{*}+\delta B\right)C$.
Hence if $\Delta K$ satisfies equation (41) then $\ensuremath{K^{*}+\Delta K\in\mathcal{K}_{P,s}\left(t^{*}+\delta t\right)}$.
Now all we need to do is to prove that $\Delta K$ is infinitesimal,
i.e. $\Delta K\rightarrow0$ if $\delta A\rightarrow0$ and $\delta B\rightarrow0$.
Taking norm on both side of equation (41) we get,

\vspace{-1.2em}

\begin{eqnarray}
\Delta K\left\Vert \mathcal{R}\right\Vert  & \leq & 2\left\Vert P\right\Vert \left(\left\Vert \delta A\right\Vert +K^{*}\left\Vert \delta B\right\Vert \left\Vert C\right\Vert \right)\nonumber \\
\Delta K & \leq & \frac{2\left\Vert P\right\Vert \left(\delta_{A}+K^{*}\delta_{B}\right)}{\left\Vert \mathcal{R}\right\Vert }\delta t
\end{eqnarray}

Since $\left\Vert \mathcal{R}\right\Vert $ is finite it is obvious
from inequality (42) that $\Delta K\rightarrow0$ as $\delta t\rightarrow0$.
This concludes the proof.

\section{Global Solution of BMI-EMP\label{sec:BMIEMP}}

\textit{Theorem 2} involves solving the following BMI-EMP optimization
problem

\noindent \textbf{OP1:}\\
\textit{minimize: }$s$\\
\textit{subject to:}\\
$\left(A_{i}-B_{i}K_{i}C\right)^{T}P+P\left(A_{i}-B_{i}K_{i}C\right)\preceq sI\quad\forall\left(A_{i},B_{i}\right)\in\mathcal{V}$
$k_{m}\leq K_{i}\leq k_{M},\quad i=1,2,\ldots,m$\\
$P\succ0$

We will first justify why \textbf{OP1} is called BMI-EMP. Consider
a matrix inequality of type

\noindent \vspace{-1.8em}

\begin{equation}
\left(A-BKC\right)^{T}P+P\left(A-BKC\right)\preceq sI
\end{equation}

Given a set of $A$, $B$, $C$, $K$ and $P$, the least possible
value of $s$ which will satisfy inequality (43) is indeed the largest
eigen value of the matrix $\left(A-BKC\right)^{T}P+P\left(A-BKC\right)$.
So if we exclude the inequalities $k_{m}\leq K_{i}\leq k_{M}$ and
$P\succ0$, then \textbf{OP1} can be equaly casted as

\noindent \vspace{-2.2em}

\[
\underset{P,\, K_{i}}{\min}\;\underset{\left(A_{i},\, B_{i}\right)\in\mathcal{V}}{\max}\lambda_{M}\left(\left(A_{i}-B_{i}K_{i}C\right)^{T}P+P\left(A_{i}-B_{i}K_{i}C\right)\right)
\]

which is basically a Largest Eigenvalue Minimization Problem or just
``EMP''. Now consider the function

\noindent \vspace{-1.8em}

\[
\Lambda\left(P\,,\, K\right)=\lambda_{M}\left(\left(A-BKC\right)^{T}P+P\left(A-BKC\right)\right)
\]

The matrix $Q\left(P,\, K\right):=\left(A-BKC\right)^{T}P+P\left(A-BKC\right)$
is \textit{Bilinear} in the sense that it is linear in $K$ if $P$
is fixed and linear in $P$ if $K$ is fixed. As $Q\in\mathbf{S}$,
the function $\Lambda\left(P\,,\, K\right)=\lambda_{M}\left(Q\left(P,\, K\right)\right)$
is convex in $K$ if $P$ is fixed and convex in $P$ if $K$ is fixed.

We are interested in the global solution of BMI-EMP as a smaller value
$s$ will ensure better convergence of the LTV/LTI system. In the
following we will provide a sketch of the work done in \cite{goh1994global,globalBMI07_parallel,globalBMI97_fujioka}
which will be just enough to design an algorithm for global solution
of BMI-EMP. However we don't provide detailed explanation of the algorithm
for which the reader may refer \cite{goh1994global}.

Before proceeding forward we would like to cast \textbf{OP1} in a
form which can be handled by a numerical solver. Observe that the
third constrain of \textbf{OP1} is a strict inequality which demands
that the Lyapunov Matrix $P$ has to be positive definite (\textit{not}
positive semi-definite). Such a strict inequality will impose numerical
challenge and hence we replace it with the non-strict inequality $P\succeq\mu_{p}I$
where, $0<\mu_{p}\ll1$. Without any loss of generality: 1) We constrain
the $\left\Vert P\right\Vert \leq1$ by imposing the constrain $P\preceq I$.
2) We normalize the RGS gain set. We get the following optimization
problem,

\noindent \textbf{OP2:}\\
\textit{minimize: }$s$\\
\textit{subject to:}\\
$A_{i}^{T}P+PA_{i}-\left(B_{i}C\right)^{T}K_{i}P-K_{i}P\left(B_{i}C\right)\preceq sI\;\forall\left(A_{i},B_{i}\right)\in\mathcal{V}$\\
$\mu_{k}\leq K_{i}\leq 1,\quad i=1,2,\ldots,m$\\
$\mu_{p}I\preceq P\preceq I$

In the above problem, $\mu_{k}=\frac{k_{m}}{k_{M}}<1$. We also redefine
$C$ as $C:=\left[\begin{array}{cccc}
k_{M} & \cdots & 0 & 0\end{array}\right]$, to neutralize the effect of normalizing RGS gain set.

We now define two vectors: $\mathbf{P}=\left[\begin{array}{cccc}
p_{1} & p_{2} & \cdots & p_{n_{p}}\end{array}\right]^{T}$ containing the $n_{p}=\frac{N\left(N+1\right)}{2}$ distinct elements
of symmetric matrix $P$ and $\mathbf{K}=\left[\begin{array}{cccc}
K_{1} & K_{2} & \cdots & K_{m}\end{array}\right]^{T}$ the $m$ normalized RGS Gains of \textbf{OP2}. We also define two
sets $\mathcal{X}_{P}$ and $\mathcal{X}_{K}$ as follows

\noindent \vspace{-2.2em}

\[
\mathcal{X}_{P}:=\left[-1,\,1\right]^{n_{p}}\:;\quad\mathcal{X}_{K}:=\left[\mu_{k},\,1\right]^{m}
\]

Note that $\mathcal{X}_{P}$ is a $n_{p}$ dimensional unit hypercube
such that%
\footnote{All the element of the Lyapunov Matrix $P$ will be in the range $\left[-1,\,1\right]$
as $P\preceq I$ according to \textbf{OP2}.%
} $\mathbf{P}\in\mathcal{X}_{P}$ and $\mathcal{X}_{K}$ is a $m$
dimensional hyper-rectangle such that $\mathbf{K}\in\mathcal{X}_{K}$.
We also define smaller hyper-rectangle's $\mathcal{Q}_{P}\subseteq\mathcal{X}_{P}$,
$\mathcal{Q}_{K}\subseteq\mathcal{X}_{K}$ and $\mathcal{Q}\subseteq\mathcal{X}_{P}\times\mathcal{X}_{K}$
as follows

\noindent \vspace{0.1em}

$\mathcal{Q}_{P}\::=\:\left[L_{P}^{1},\, U_{P}^{1}\right]\times\left[L_{P}^{2},\, U_{P}^{2}\right]\times\ldots\times\left[L_{P}^{n_{p}},\, U_{P}^{n_{p}}\right]$\\

$\mathcal{Q}_{K}\::=\:\left[L_{K}^{1},\, U_{K}^{1}\right]\times\left[L_{K}^{2},\, U_{K}^{2}\right]\times\ldots\times\left[L_{K}^{m},\, U_{K}^{m}\right]$\\

$\mathcal{Q}\::=\:\mathcal{Q}_{P}\times\mathcal{Q}_{K}$

Obviously $-1\leq L_{P}^{i}\leq U_{P}^{i}\leq1$ and $\mu_{k}\leq L_{K}^{i}\leq U_{K}^{i}\leq1$.
Consider the following ``constrained version'' of \textbf{OP2}
where $\ensuremath{\left(\mathbf{P}\,,\,\mathbf{K}\right)}$ is only
defined in the small hyper-rectangle $\mathcal{Q}$.

\noindent \textbf{OP3:}\\
\textit{minimize: }$s$\\
\textit{subject to:}\\
$A_{i}^{T}P+PA_{i}-\left(B_{i}C\right)^{T}K_{i}P-K_{i}P\left(B_{i}C\right)\preceq sI\;\forall\left(A_{i},B_{i}\right)\in\mathcal{V}$\\
$\mu_{k}\leq K_{i}\leq 1,\quad i=1,2,\ldots,m$\\
$\mu_{p}I\preceq P\preceq I$\\
$\left(\mathbf{P}\,,\,\mathbf{K}\right)\in\mathcal{Q}$

The input to \textbf{OP3} is the set $\mathcal{Q}$ and its output
is $s^{*}\left(\mathcal{Q}\right)$ which is a function of $\mathcal{Q}$.
We want to bound $s^{*}\left(\mathcal{Q}\right)$ by an upper and
a lower bound as follows:

\noindent \vspace{-2.0em}

\begin{equation}
\Phi_{L}\left(\mathcal{Q}\right)\leq s^{*}\left(\mathcal{Q}\right)\leq\Phi_{U}\left(\mathcal{Q}\right)
\end{equation}

\noindent The \textit{convex relaxation} technique introduced in
\cite{globalBMI97_fujioka} (also discussed in \cite{globalBMI07_parallel},
equation (9)) has been used to get the lower bound $\Phi_{L}\left(\mathcal{Q}\right)$.
We replace the nonlinearity $K_{i}P$ with a new matrix $W_{i}$.
As $W_{i}$ is a symmetric matrix matrix of order $N$ we can represent
it by a vector $\mathbf{W}_{i}=\left[\begin{array}{cccc}
w_{1i} & w_{2i} & \cdots & w_{n_{p}i}\end{array}\right]^{T}$ containing the $n_{p}$ distinct elements of $W_{i}$. We now define
a matrix $\mathbf{W}=\left[\begin{array}{cccc}
\mathbf{W}_{1}^{T} & \mathbf{W}_{2}^{T} & \cdots & \mathbf{W}_{m}^{T}\end{array}\right]^{T}$. If \textbf{OP3} is expanded in terms of $p_{j}$ and $K_{i}$ then
the element $w_{ji}$ of $\mathbf{W}$ is constrained by the equality
$w_{ji}=K_{i}p_{j}$. Rather than imposing this equality constrain
we let $w_{ji}$ to be a free variable which can take any value in
the set $\mathcal{W}\left(\mathcal{Q}\right)$ defined as

\noindent \vspace{-2.0em}

\[
\mathcal{W}\left(\mathcal{Q}\right):=\left\{ \mathbf{W}\left|\begin{array}{c}
\ensuremath{\left(\mathbf{P}\,,\,\mathbf{K}\right)\in\mathcal{Q}}\\
w_{ji}\geq L_{K}^{i}p_{j}+L_{P}^{j}K_{i}-L_{K}^{i}L_{P}^{j}\\
w_{ji}\geq U_{K}^{i}p_{j}+U_{P}^{j}K_{i}-U_{K}^{i}U_{P}^{j}\\
w_{ji}\leq U_{K}^{i}p_{j}+L_{P}^{j}K_{i}-U_{K}^{i}L_{P}^{j}\\
w_{ji}\leq L_{K}^{i}p_{j}+U_{P}^{j}K_{i}-L_{K}^{i}U_{P}^{j}
\end{array}\right.\right\} 
\]

By performing the convex relaxation stated above we get the following
optimization problem.

\noindent \textbf{OP4:}\\
\textit{minimize: }$s$\\
\textit{subject to:}\\
$A_{i}^{T}P+PA_{i}-\left(B_{i}C\right)^{T}W_{i}-W_{i}\left(B_{i}C\right)\preceq sI\quad\forall\left(A_{i},B_{i}\right)\in\mathcal{V}$\\
$\mu_{p}I\preceq P\preceq I$\\
$\mu_{p}K_{i}I\mu_{k}\preceq W_{i}\preceq K_{i}I,\quad i=1,2,\ldots,m$\\
$\left(\mathbf{P}\,,\,\mathbf{K}\right)\in\mathcal{Q}$\\
$\mathbf{W}\in\mathcal{W}\left(\mathcal{Q}\right)$

The input to \textbf{OP4} is the set $\mathcal{Q}$ and its output
is $\Phi_{L}\left(\mathcal{Q}\right)$. As \textbf{OP4} is a relaxed
version of \textbf{OP3}, $\Phi_{L}\left(\mathcal{Q}\right)\leq s^{*}\left(\mathcal{Q}\right)$.
Note that \textbf{OP4} is a convex optimization problem, more specifically
a \textit{Semi-Definite Program(SDP)} which can be solved by numerical
solvers like CVX\cite{cvx}.

\begin{algorithm}[t]
1. Set $\mathbf{K}^{\left(0\right)}=$Centroid of $\mathcal{Q}_{K}$.$\left(s^{\left(0\right)},\mathbf{P}^{\left(0\right)}\right)=\mathbf{OP3}\left(\mathcal{Q},\mathbf{K}^{\left(0\right)}\right)$.

2. $if\quad\left(s^{\left(0\right)}=\infty\right)$

3.$\qquad$$return$ $\infty$.

4. $else$

5.$\qquad$Set $\delta>0$ and $k=0$.

6.$\qquad$$do\:\{$

7.$\qquad$$\qquad$$\quad\:$$\left(s^{\left(k+1\right)}\,,\,\mathbf{P}^{\left(k+1\right)}\right)=\mathbf{OP3}\left(\mathcal{Q}\,,\,\mathbf{K}^{\left(k\right)}\right)$

8.$\qquad$$\qquad$$\quad\:$$\left(s^{\left(k+1\right)}\,,\,\mathbf{K}^{\left(k+1\right)}\right)=\mathbf{OP3}\left(\mathcal{Q}\,,\,\mathbf{P}^{\left(k+1\right)}\right)$

9.$\qquad$$\qquad$$\quad\:$$k=k+1$.

10.$\qquad$$\}\: while\:\left(s^{\left(k-1\right)}-s^{\left(k\right)}<\delta\left|s^{\left(k\right)}\right|\right)$

11.$\qquad$$return$ $s^{\left(k\right)}$.

12. $end$

\caption{Alternating SDP Method}
\end{algorithm}

We now concentrate on defining the upper bound $\Phi_{U}\left(\mathcal{Q}\right)$.
Any local minima of \textbf{OP3} can indeed be the upper bound $\Phi_{U}\left(\mathcal{Q}\right)$.
We use the \textit{Alternating SDP Method} discussed in \cite{globalBMI97_fujioka,globalBMI97_Fukuda}
. Alternating SDP method relies on the Bi-Convex nature of \textbf{OP3},
i.e. \textbf{OP3} becomes a convex problem (more specifically SDP)
in $\mathbf{P}$ with $\mathbf{K}$ fixed or a convex problem in $\mathbf{K}$
with $\mathbf{P}$ fixed. Alternating SDP Method is summarized in
\textit{Algorithm 1}. We represent by $\mathbf{OP3}\left(\mathcal{Q}\,,\,\mathbf{K}'\right)$
the optimization problem obtained by fixing $\mathbf{K}=\mathbf{K}^{'}$
in \textbf{OP3} and $\mathbf{OP3}\left(\mathcal{Q}\,,\,\mathbf{P}'\right)$
the optimization problem obtained by fixing $\mathbf{P}=\mathbf{P}^{'}$
in \textbf{OP3}. The input to $\mathbf{OP3}\left(\mathcal{Q}\,,\,\mathbf{K}'\right)$
is the small hyper-rectangle $\mathcal{Q}$ and the fixed RGS gain
$\mathbf{K}^{'}$ while the input to $\mathbf{OP3}\left(\mathcal{Q}\,,\,\mathbf{P}'\right)$
is the small hyper-rectangle $\mathcal{Q}$ and the fixed Lyapunov
Matrix $\mathbf{P}^{'}$. The outputs of $\mathbf{OP3}\left(\mathcal{Q}\,,\,\mathbf{K}'\right)$
are the optimized $s$ and $\mathbf{P}$ while the outputs of $\mathbf{OP3}\left(\mathcal{Q}\,,\,\mathbf{P}'\right)$
are the optimized $s$ and $\mathbf{K}$.

Under the various definations introduced above the \textit{Branch and Bound Algorithm}
\cite{goh1994global} calculate the global minima of \textbf{OP2}
to an absolute accuracy of $\epsilon>0$ within finite time. The psuedocode
of \textit{Branch and Bound} method is given in \textit{Algorithm 2}.

\begin{algorithm}[t]
1. Set $\epsilon>0$ and $k=0$.

2. Set $\mathcal{Q}_{0}=\mathcal{X}_{P}\times\mathcal{X}_{K}$ and
$\mathcal{G}_{0}=\left\{ \mathcal{Q}_{0}\right\} $.

3. Set $L_{0}=\Phi_{L}\left(\mathcal{Q}_{0}\right)$ and $U_{0}=\Phi_{U}\left(\mathcal{Q}_{0}\right)$.

4. $while\quad\left(U_{k}-L_{k}<\epsilon\right)$

5.$\qquad$Select $\bar{\mathcal{Q}}$ from $\mathcal{G}_{k}$ such
that $L_{k}=\Phi_{L}\left(\mathcal{\bar{Q}}\right)$.

6.$\qquad$Set $\mathcal{G}_{k+1}=\mathcal{G}_{k}-\left\{ \mathcal{\bar{Q}}\right\} $.

7.$\qquad$Split $\bar{\mathcal{Q}}$ along its longest egde into
$\bar{\mathcal{Q}}_{1}$ and $\bar{\mathcal{Q}}_{2}$.

8.$\qquad$$for\quad\left(i=1,\:2\right)$

9.$\qquad\quad$$if\quad\left(\Phi_{L}\left(\bar{\mathcal{Q}}_{i}\right)\leq U_{k}\right)$

10.$\qquad\qquad$Compute $\Phi_{U}\left(\bar{\mathcal{Q}}_{i}\right)$.

11.$\qquad\qquad$Set $\mathcal{G}_{k+1}=\mathcal{G}_{k}\bigcup\left\{ \bar{\mathcal{Q}}_{i}\right\} $.

12.$\qquad\:\,$$end$

13.$\quad\;$$end$

14.$\quad\;$Set $U_{k+1}=\underset{\mathcal{Q}\in\mathcal{G}_{k+1}}{min}\:\Phi_{U}\left(\mathcal{Q}\right)$.

15.$\quad\;$$Pruning\,:\;\mathcal{G}_{k+1}=\mathcal{G}_{k+1}-\left\{ \mathcal{Q}\,:\:\Phi_{L}\left(\mathcal{Q}\right)>U_{k+1}\right\} $.

16.$\quad\;$Set $L_{k+1}=\underset{\mathcal{Q}\in\mathcal{G}_{k+1}}{min}\:\Phi_{L}\left(\mathcal{Q}\right)$.

17.$\quad\;$Set $k=k+1$.

18. $end$

\caption{Branch and Bound}
\end{algorithm}

\section{Control Theoretic Definitions and Concepts\label{sec:Notions}}

The control theoretic analysis presented in this paper relies on definitions
and theorems from three broad areas. This concepts are standard and
can be easily found in books like \cite{kreyszig1989introductory,boyd2009convex}.

\subsection{Real Analysis}

\textit{Defintion 1 (Cartesian Product of Sets): }Cartesian Products
of two sets $A$ and $B$, denoted by $A\times B$, is defined as
the set of all ordered pairs $\left(a,\, b\right)$ where $a\in A$
and $b\in B$.

\textit{Notion 1 (Connected Set): }A set is said to be connected
if for any two points on that set, there exist \textit{at-least}
one path joining those two points which also lies on the set. \textit{Note that this is not the "general definition" of connected set.}

\textit{Notion 2 (Compact Set): }In Eucledian Space $\mathbb{R}^{n}$,
a closed and bounded set is called a compact set. Also a compact set
in $\mathbb{R}^{n}$ is always closed and bounded. \textit{Note that this is not the "general definition" of compact set.}

\subsection{Linear Algebra }

\textit{Definition 2 (Eucledian Norm of a Vector): }For a vector
$x\in\mathbb{R}^{n}$, eucledian norm $\left\Vert x\right\Vert $
is defined as $\left\Vert x\right\Vert :=\sqrt{x^{T}x}$. Throughout
the entire paper ``norm'' means ``eucledian norm'' unless mentioned
otherwise.

\textit{Definition 3 (Induced Spectral Norm of a Matrix): }For a
matrix $A\in\mathbb{R}^{m\times n}$, Induced Spectral Norm $\left\Vert A\right\Vert $
is defined as

\[
\left\Vert A\right\Vert :=\underset{x\in\mathbb{R}^{n}-\left\{ 0\right\} }{\sup}\frac{\left\Vert Ax\right\Vert }{\left\Vert x\right\Vert }
\]

It can equally be defined as

\noindent \vspace{-1.4em}

\[
\left\Vert A\right\Vert :=\underset{\left\Vert x\right\Vert =1}{\sup}\left\Vert Ax\right\Vert 
\]

\vspace{-0.6em}

\textit{Definition 4 (Positive (Negative) Definite (Semi-Definite) Matrix): }A
square matrix $A\in\mathbb{R}^{n\times n}$ is said to be positive
definite if $x^{T}Ax>0,\:\forall\, x\in\mathbb{R}^{n}-\left\{ 0\right\} $
and is said to be positive semi-definite if $x^{T}Ax\geq0,\:\forall\, x\in\mathbb{R}^{n}$.
A matrix $A$ is said to be negative definite is $-A$ is positive
definite and is said to be negative semi-definite is $-A$ is positive
semi-definite.

\textit{Linear Algebra Theorems:}\begin{enumerate}
\item Properties of norms:
\begin{enumerate}
\item For a vector $x\in \mathbb{R}^{n}$, if $\left\Vert x \right\Vert=0$ then $x=0$. Also if $\left\Vert x \right\Vert\rightarrow 0$ then $x\rightarrow 0$.
\item For any two matrix $A,\,B\in\mathbb{R}^{n \times m}$, $\left\Vert A+B\right\Vert \leq\left\Vert A\right\Vert +\left\Vert B\right\Vert $.
\item For any matrix $A\in\mathbb{R}^{n \times m}$ and a scalar $\alpha$, $\left\Vert \alpha A\right\Vert=\alpha \left\Vert A\right\Vert$.
\item For any two matrix $A,\,B\in\mathbb{R}^{n \times m}$, $\left\Vert AB\right\Vert \leq\left\Vert A\right\Vert \left\Vert B\right\Vert $.
\item For any matrix $A\in\mathbb{R}^{n \times m}$ and a vector $x\in\mathbb{R}^{m}$, $\left\Vert Ax\right\Vert \leq\left\Vert A\right\Vert \left\Vert x\right\Vert $.
\end{enumerate}
\item Eigenvalues of symmetric positive (negative) definite matrix are positive (negative).
\item For any symmetric matrix $A\in\mathbb{R}^{n \times n}$ and a vector $x\in\mathbb{R}^{n}$,
\[ \lambda_{m}\left(A\right)\left\Vert x\right\Vert ^{2}\leq x^{T}Ax\leq\lambda_{M}\left(A\right)\left\Vert x\right\Vert ^{2} \]
\item For any symmetric positive definite (semi-definite) matrix $A$, $\left\Vert A\right\Vert =\lambda_{M}\left(A\right)$. 
\item For any square matrix $A$, $x^{T}Ax\leq\left\Vert A\right\Vert \left\Vert x\right\Vert ^{2}$
\end{enumerate}

\subsection{Convex Analysis}

\textit{Definition 5 (Convex Set): }A set $\mathcal{A}$ is said
to be convex if for any $x\,,\, y\in\mathcal{A}$, $\theta x+\left(1-\theta\right)y\in\mathcal{A}\,,\:\forall\,\theta\in\left[0\,,\,1\right]$.

\textit{Definition 6 (Convex Function): }Let $\mathcal{A}\in\mathbb{R}^{n}$
be a convex set. A function $f\,:\,\mathcal{A}\rightarrow\mathbb{R}$
is convex if

\vspace{-0.6em}

\[
f\left(\theta x+\left(1-\theta\right)y\right)\leq\theta f\left(x\right)+\left(1-\theta\right)f\left(y\right)
\]

\noindent for all $x\,,\, y\in\mathcal{A}$ and for all $0\leq\theta\leq1$.

\textit{Definition 7 (Convex Combination and Convex Hull): }Given
a set $\mathcal{A}=\left\{ \begin{array}{cccc}
a_{1} & a_{2} & \cdots & a_{n}\end{array}\right\} $, its convex combination are those elements which can be expressed
as

\vspace{-0.6em}

\[
\theta_{1}a_{1}+\theta_{2}a_{2}+\cdots+\theta_{n}a_{n}
\]

where $\theta_{1}+\theta_{2}+\cdots+\theta_{n}=1$ and $\theta_{i}\geq0,\: i=1,\,2,\ldots,n$.

The set of all convex combination of set $\mathcal{A}$ is called
the convex hull of $\mathcal{A}$. Convex Hull of set $\mathcal{A}$
is indeed the smallest convex set containing $\mathcal{A}$.

\textit{Definition 8 (Semi-Definite Program): }A semi-definite program
or SDP is an optimization problem in variable $x=\begin{bmatrix}x_{1} & x_{2} & \cdots & x_{n}\end{bmatrix}^{T}\in\mathbb{R}^{n}$
with the generic structure

\[
\begin{array}{l}
minimize:\quad c^{T}x\\
subject\; to:\\
F_{0}+x_{1}F_{1}+x_{2}F_{2}+\cdots+x_{n}F_{n}\preceq0\\
Ax=b
\end{array}
\]

where $F_{0},\, F_{1},\, F_{2},\ldots,\, F_{n}\in\mathbf{S}$, $c\in\mathbb{R}^{n}$,
$A\in\mathbb{R}^{p\times n}$ and $b\in\mathbb{R}^{p}$.

\section*{Acknowledgment}

We would like to thank Prof. Radha Krishna Ganti for providing invaluable
help in areas related to convex optimization.

\bibliographystyle{IEEEtran}
\bibliography{file}

\begin{thebibliography}{10}
\providecommand{\url}[1]{#1}
\csname url@samestyle\endcsname
\providecommand{\newblock}{\relax}
\providecommand{\bibinfo}[2]{#2}
\providecommand{\BIBentrySTDinterwordspacing}{\spaceskip=0pt\relax}
\providecommand{\BIBentryALTinterwordstretchfactor}{4}
\providecommand{\BIBentryALTinterwordspacing}{\spaceskip=\fontdimen2\font plus
\BIBentryALTinterwordstretchfactor\fontdimen3\font minus
  \fontdimen4\font\relax}
\providecommand{\BIBforeignlanguage}[2]{{%
\expandafter\ifx\csname l@#1\endcsname\relax
\typeout{** WARNING: IEEEtran.bst: No hyphenation pattern has been}%
\typeout{** loaded for the language `#1'. Using the pattern for}%
\typeout{** the default language instead.}%
\else
\language=\csname l@#1\endcsname
\fi
#2}}
\providecommand{\BIBdecl}{\relax}
\BIBdecl

\bibitem{CHUA1971}
L.~Chua, ``Memristor-the missing circuit element,'' \emph{IEEE Trans. Circuit
  Theory}, vol.~18, no.~5, pp. 507--519, 1971.

\bibitem{STRUKOV2008}
D.~B. Strukov, G.~S. Snider, D.~R. Stewart, and R.~S. Williams, ``The missing
  memristor found,'' \emph{Nature}, vol. 453, no. 7191, pp. 80--83, 2008.

\bibitem{MEMORY2011}
Y.~Ho, G.~M. Huang, and P.~Li, ``Dynamical properties and design analysis for
  nonvolatile memristor memories,'' \emph{IEEE Trans. Circuits Syst. I, Reg.
  Papers}, vol.~58, no.~4, pp. 724--736, 2011.

\bibitem{JOGELKAR2009}
Y.~N. Joglekar and S.~J. Wolf, ``The elusive memristor: properties of basic
  electrical circuits,'' \emph{Eur. J. Phys.}, vol.~30, no.~4, pp. 661--675,
  2009.

\bibitem{BCM2012}
F.~Corinto and A.~Ascoli, ``A boundary condition-based approach to the modeling
  of memristor nanostructures,'' \emph{IEEE Trans. Circuits Syst. I, Reg.
  Papers}, vol.~60, pp. 2713--2726, 2012.

\bibitem{TEAM2013}
S.~Kvatinsky, E.~G. Friedman, A.~Kolodny, and U.~C. Weiser, ``Team: Threshold
  adaptive memristor model,'' \emph{IEEE Trans. Circuits Syst. I, Reg. Papers},
  vol.~60, no.~1, pp. 211--221, 2013.

\bibitem{SynapseMemristor}
S.~H. Jo, T.~Chang, I.~Ebong, B.~B. Bhadviya, P.~Mazumder, and W.~Lu,
  ``Nanoscale memristor device as synapse in neuromorphic systems,'' \emph{Nano
  letters}, vol.~10, no.~4, pp. 1297--1301, 2010.

\bibitem{memneural}
Y.~V. Pershin and M.~Di~Ventra, ``Experimental demonstration of associative
  memory with memristive neural networks,'' \emph{Neural Networks}, vol.~23,
  no.~7, pp. 881--886, 2010.

\bibitem{CharacterRecognition}
A.~Sheri, H.~Hwang, M.~Jeon, and B.~Lee, ``Neuromorphic character recognition
  system with two pcmo-memristors as a synapse,'' \emph{IEEE Trans. Ind.
  Electron.}, vol.~61, pp. 2933--2941, 2014.

\bibitem{amoebamemristivelearning}
Y.~V. Pershin, S.~La~Fontaine, and M.~Di~Ventra, ``Memristive model of amoeba
  learning,'' \emph{Phys. Rev. E}, vol.~80, no.~2, p. 021926, 2009.

\bibitem{harmonic_generation}
G.~Z. Cohen, Y.~V. Pershin, and M.~Di~Ventra, ``Second and higher harmonics
  generation with memristive systems,'' \emph{Appl. Phys. Lett.}, vol. 100,
  no.~13, p. 133109, 2012.

\bibitem{progAnalog1}
Y.~V. Pershin and M.~Di~Ventra, ``Practical approach to programmable analog
  circuits with memristors,'' \emph{IEEE Trans. Circuits Syst. I, Reg. Papers},
  vol.~57, no.~8, pp. 1857--1864, 2010.

\bibitem{progAnalog2}
S.~Shin, K.~Kim, and S.~Kang, ``Memristor applications for programmable analog
  ics,'' \emph{IEEE Trans. Nanotechnol.}, vol.~10, no.~2, pp. 266--274, 2011.

\bibitem{jeltsema2012port}
D.~Jeltsema and A.~Doria-Cerezo, ``Port-hamiltonian formulation of systems with
  memory,'' \emph{Proc. IEEE}, vol. 100, pp. 1928--1937, 2012.

\bibitem{jeltsema2010memristive}
D.~Jeltsema and A.~J. van~der Schaft, ``Memristive port-hamiltonian systems,''
  \emph{Math. Comput. Model. Dyn. Syst.}, vol.~16, pp. 75--93, 2010.

\bibitem{stabilityMemNeural}
J.~Hu and J.~Wang, ``Global uniform asymptotic stability of memristor-based
  recurrent neural networks with time delays,'' in \emph{Proc. Int. Joint Conf.
  Neural Netw.}, 2010, pp. 1--8.

\bibitem{DELGADO2010}
A.~Delgado, ``The memristor as controller,'' in \emph{IEEE Nanotechnology
  Materials and Devices Conference}, 2010, pp. 376--379.

\bibitem{MEMIDAPBC2012}
J.~S. A.~Doria-Cerezo, L. van der~Heijden, ``Memristive port-hamiltonian
  control: Path-dependent damping injection in control of mechanical systems,''
  in \emph{Proc. 4th IFAC Workshop on Lagrangian and Hamiltonian Methods for
  Non Linear Control}, 2012, pp. 167--172.

\bibitem{cdc_memristor}
R.~Pasumarthy, G.~Saha, F.~Kazi, and N.~Singh, ``Energy and power based
  perspective of memristive controllers,'' in \emph{Proc. IEEE Conf. Decis.
  Contr.}, 2013, pp. 642--647.

\bibitem{shapiro2011feedback}
J.~J. Gorman and B.~Shapiro, \emph{Feedback control of mems to atoms}.\hskip
  1em plus 0.5em minus 0.4em\relax Springer, 2011.

\bibitem{tunableresistor}
E.~Ozalevli and P.~E. Hasler, ``Tunable highly linear floating-gate cmos
  resistor using common-mode linearization technique,'' \emph{IEEE Trans.
  Circuits Syst. I, Reg. Papers}, vol.~55, pp. 999--1010, 2008.

\bibitem{KANG1976}
L.~O. Chua and S.~M. Kang, ``Memristive devices and systems,'' \emph{Proc.
  IEEE}, vol.~64, no.~2, pp. 209--223, 1976.

\bibitem{memgain}
T.~A. Wey and W.~D. Jemison, ``Variable gain amplifier circuit using titanium
  dioxide memristors,'' \emph{IET Circuits, Devices \& Syst.}, vol.~5, no.~1,
  pp. 59--65, 2011.

\bibitem{haykin2009communication}
S.~Haykin, \emph{Communication systems}, 5th~ed.\hskip 1em plus 0.5em minus
  0.4em\relax Wiley Publishing, 2009.

\bibitem{boyd1994linear}
S.~P. Boyd, \emph{Linear matrix inequalities in system and control
  theory}.\hskip 1em plus 0.5em minus 0.4em\relax SIAM, 1994.

\bibitem{ben2001lectures}
A.~Ben-Tal and A.~Nemirovski, \emph{Lectures on modern convex optimization:
  analysis, algorithms, and engineering applications}.\hskip 1em plus 0.5em
  minus 0.4em\relax SIAM, 2001.

\bibitem{rosenbrook}
H.~H. Rosenbrook, ``The stability of linear time-dependent control systems,''
  \emph{Int. J. Electron. \& Contr.}, vol.~15, pp. 73--80, 1963.

\bibitem{solo}
V.~Solo, ``On the stability of slowly time-varying linear systems,''
  \emph{Math. of Contr., Signals \& Syst.}, vol.~7, pp. 331--350, 1994.

\bibitem{qhull}
\BIBentryALTinterwordspacing
C.~B. Barber, D.~P. Dobkin, and H.~Huhdanpaa, ``The quickhull algorithm for
  convex hulls,'' \emph{ACM Trans. on Mathematical Software}, vol.~22, no.~4,
  pp. 469--483, 1996. [Online]. Available: \url{http://www.qhull.org}
\BIBentrySTDinterwordspacing

\bibitem{localBMI99_pathfollowing}
A.~Hassibi, J.~How, and S.~Boyd, ``A path-following method for solving bmi
  problems in control,'' in \emph{Proc. Amer. Contr. Conf.}, 1999, pp.
  1385--1389.

\bibitem{localBMI98_ILMI}
Y.-Y. Cao, J.~Lam, and Y.-X. Sun, ``Static output feedback stabilization: an
  ilmi approach,'' \emph{Automatica}, vol.~34, no.~12, pp. 1641--1645, 1998.

\bibitem{goh1994global}
K.-C. Goh, M.~Safonov, and G.~Papavassilopoulos, ``A global optimization
  approach for the bmi problem,'' in \emph{Proc. IEEE Conf. Decis. Contr},
  vol.~3, 1994, pp. 2009--2014.

\bibitem{globalBMI97_fujioka}
K.~H. H.~Fujioka, ``Bounds for the bmi eingenvalue problem - a good lower bound
  and a cheap upper bound,'' \emph{Trans. Society of Instr. \& Contr. Eng.},
  vol.~33, pp. 616--621, 1997.

\bibitem{globalBMI07_parallel}
M.~Michihiro~Kawanishi and Y.~Shibata, ``Bmi global optimization using parallel
  branch and bound method with a novel branching method,'' in \emph{Proc. Amer.
  Contr. Conf.}, 2007, pp. 1664--1669.

\bibitem{ogata}
K.~Ogata, \emph{Modern Control Engineering}, 5th~ed.\hskip 1em plus 0.5em minus
  0.4em\relax Prentice Hall, 2010.

\bibitem{boyd2009convex}
S.~Boyd and L.~Vandenberghe, \emph{Convex optimization}.\hskip 1em plus 0.5em
  minus 0.4em\relax Cambridge university press, 2009.

\bibitem{rockafellar1997convex}
R.~T. Rockafellar, \emph{Convex analysis}.\hskip 1em plus 0.5em minus
  0.4em\relax Princeton university press, 1997.

\bibitem{sparse98_Hassibi}
A.~Hassibi, J.~How, and S.~Boyd, ``Low-authority controller design via convex
  optimization,'' in \emph{Proc. IEEE Conf. Decis. Contr.}, 1998, pp.
  1385--1389.

\bibitem{barkana2013simple}
I.~Barkana, ``Simple adaptive control-a stable direct model reference adaptive
  control methodology-brief survey,'' \emph{Int. J. Adapt. Contr. Signal
  Proc.}, 2013.

\bibitem{cvx}
I.~CVX~Research, ``{CVX}: Matlab software for disciplined convex programming,
  version 2.0,'' \url{http://cvxr.com/cvx}, Aug. 2012.

\bibitem{globalBMI97_Fukuda}
M.~K. Mituhiro~Fukuda, ``Branch-and-cut algorithms for the bilinear matrix
  inequality eigenvalue problem,'' \emph{Comput. Optim. Appl}, vol.~19, p.
  2001, 1999.

\bibitem{kreyszig1989introductory}
E.~Kreyszig, \emph{Introductory functional analysis with applications}.\hskip
  1em plus 0.5em minus 0.4em\relax Wiley, 1989.

\end{thebibliography}

\vspace{-3.0em}
\begin{IEEEbiography}[{\includegraphics[clip,width=1in,height=1.25in]{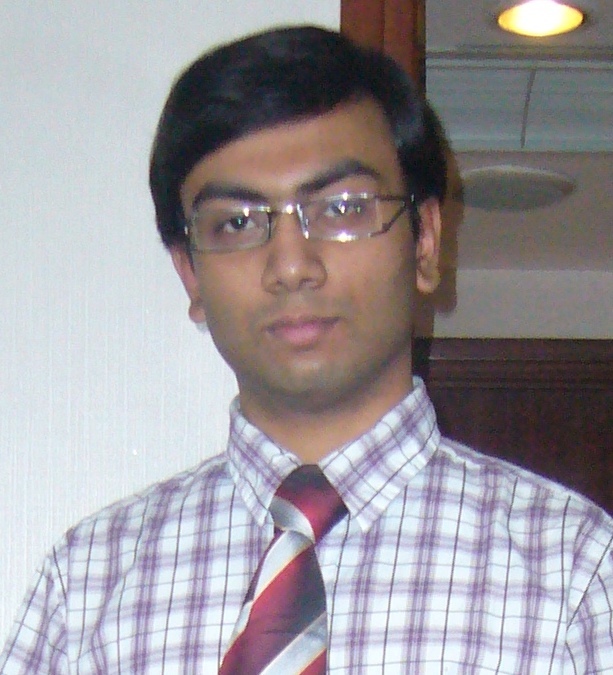}}]
{Gourav Saha} received the B.E. degree from Anna University, Chennai,
India, in 2012. Since 2013, he has been working towards his M.S. (by
research) degree in the Department of Electrical Engineering at Indian
Institute of Technology, Madras.

\noindent His current research interests include application of memristor
and memristive systems, control of MEMS and modelling and control
of cyber-physical systems with specific focus on cloud computing. 
\end{IEEEbiography}

\vspace{-3.0em}
\begin{IEEEbiography}[{\includegraphics[clip,width=1in,height=1.25in]{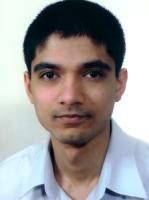}}]
{Ramkrishna Pasumarthy} obtained his PhD in Systems and Control
from the University of Twente, The Netherlands in 2006. He is currently
an Assistant Professor at the Department of Electrical Engineering,
IIT Madras, India.

His research interests lie in the area of modeling and control of
physical systems, infinite dimensional systems and control of computing
systems, with specific focus on cloud computing systems.
\end{IEEEbiography}

\vspace{-3.0em}
\begin{IEEEbiography}[{\includegraphics[clip,width=1in,height=1.25in]{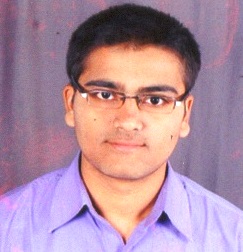}}]
{Prathamesh Khatavkar} received the B.Tech degree from University
of Pune, India, in 2012. He is currently pursuing his M.S. (by research)
degree in the Department of Electrical Engineering at Indian Institute
of Technology, Madras.

\noindent His current research interests includes Analog/RF IC design.\end{IEEEbiography}

\end{document}